\documentclass[10pt,a4paper,draft]{article}
\usepackage{bbm}
\usepackage{pifont}
\usepackage{bbding}
\usepackage{amsmath, esint}
\usepackage{mathrsfs}
\usepackage{amsfonts}
\usepackage{amssymb}
\usepackage{amsfonts,amssymb,amsmath,indentfirst,cases,amsthm}
\usepackage{epsfig}
\usepackage[numbers,sort&compress]{natbib}
\def\ess~sup{\mathop{\rm ess~sup}}

\numberwithin{equation}{section}

\newenvironment{key words}{\emph{\texttt{Keywords}}\mbox{  }}{ }
\newtheorem{theorem}{Theorem}[section]
\newtheorem{lemma}[theorem]{Lemma}

\newtheorem{definition}[theorem]{Definition}

\theoremstyle{remark}

\theoremstyle{plain}

\makeatletter

\newcommand{\Rmnum}[1]{\expandafter\@slowromancap\romannumeral #1@}
\makeatother
\begin{document}

\title{\textbf{A Nash type result for Divergence Parabolic Equation related to H\"{o}rmander's vector fields}
\thanks{This work is supported by the National Natural Science
Foundation of China (No.11271299).}}
\author{Lingling Hou ,  Pengcheng Niu \thanks{Corresponding
author's E-mail: pengchengniu@nwpu.edu.cn(P. Niu)}, \\  % Qianxia Guo, Huijv Wang%
\small{Department of Applied Mathematics, Northwestern Polytechnical
University,}\\ \small{ Xi'an, Shaanxi, 710129, P. R. China}}
\date{} \maketitle

% ----------------------------------------------------------------

\maketitle {\bf Abstract.}\
In this paper we consider the divergence parabolic equation with bounded and measurable coefficients related to H\"{o}rmander's vector fields and establish a Nash type result, i.e., the local H\"{o}lder regularity for weak solutions. After deriving the parabolic Sobolev inequality, (1,1) type Poincar\'{e} inequality of H\"{o}rmander's vector fields and a De Giorgi type Lemma, the H\"{o}lder regularity of weak solutions to the equation is proved based on the estimates of oscillations of solutions and the isomorphism between parabolic Campanato space and parabolic H\"{o}lder space. As a consequence, we give the Harnack inequality of weak solutions by showing an extension property of positivity for functions in the De Giorgi class.


\def\Xint#1{\mathchoice
    {\XXint\displaystyle\textstyle{#1}}%
    {\XXint\textstyle\scriptstyle{#1}}%
    {\XXint\scriptstyle\scriptscriptstyle{#1}}%
    {\XXint\scriptscriptstyle\scriptscriptstyle{#1}}%
    \!\int}
\def\XXint#1#2#3{{\setbox0=\hbox{$#1{#2#3}{\int}$}
    \vcenter{\hbox{$#2#3$}}\kern-.5\wd0}}
\def\dashint{\Xint-}

%---------------------------------------------------------------------------------
\section{\textbf{Introduction}\label{Section 1}}
%---------------------------------------------------------------------------------

Schauder theory for the solutions to linear elliptic and parabolic equations with $C^\alpha$ coefficients or VMO coefficients has been completed. De Girogi has followed the local H\"{o}lder continuity for the solutions to the divergence elliptic equation with bounded and measurable coefficients
\[-\sum\limits_{i,j = 1}^n {D_i}\left({a^{ij}}(x){D_j}u \right)=0, x\in{\mathbb{R}^n}\]
and given the a priori estimate of H\"{o}lder norm (see \cite{DG}). Nash in \cite{N} derived independently the similar result for the solutions to the parabolic equation with a different approach from \cite{DG}. Hereafter Moser in \cite{M} developed a new method (nowadays it is called the Moser iteration method) and proved again results above-mentioned to elliptic and parabolic equations. These important works break a new path for the study of regularity for weak solutions to partial differential equations.

In \cite{FS} Fabes and Stroock proved the Harnack inequality for linear parabolic equations by going back to Nash's original technique in \cite{N}. A very interesting approach has been raised by De Benedetto (see \cite{DB}) for proving a Harnack inequality of functions belonging to parabolic De Giorgi classes. The approach was used to derive the H\"{o}lder continuity of solutions to linear second order parabolic equations with bounded and measurable coefficients (see \cite{LSU}). Giusti \cite{G} applied the approach to give a proof of the Harnack inequality in the elliptic setting.

Square sum operators constructed by vector fields satisfying the finite rank condition were introduced by H\"{o}rmander (see \cite{H1}), who deduced that such operators are hypoelliptic. Many authors carried on researches to such operators and acquired numerous important results (\cite{Fo1},\cite{Fo2},\cite{H2},\cite{K},\cite{RS}). Nagel, Stein and Wainger (\cite{NSW}) concluded the deep properties of balls and metrics defined by vector fields of this type. Many other authors obtained very appreciable results related to H\"{o}rmander's square sum operators, for instance fundamental solutions (\cite{S}), the Poincar\'{e} inequality (\cite{J}), potential estimates (\cite{B}), also see \cite{FS},\cite{JS}, etc.. All these motivate the study to degenerate elliptic and parabolic equations formed from H\"{o}rmander's vector fields. Schauder estimates to degenerate elliptic and parabolic equations related to noncommutative vector fields have been handled in \cite{CH},\cite{GL},\cite{XZ}, etc. Bramanti and Brandolini in \cite{B2} investigated Schauder estimates to the following H\"{o}rmander type nondivergence parabolic operator
\[H=\partial_t-\sum\limits_{i,j = 1}^q {a_{ij}}(t,x){X_i}{X_j}-\sum\limits_{i= 1}^q {b_{i}}(t,x){X_i}-c(t,x),\]
where coefficients ${a_{ij}}(t,x),{b_{i}}(t,x)$ and $ c(t,x) $ are $C^\alpha$.

In this paper, we are concerned with the divergence parabolic equation with bounded and measurable coefficients related to H\"{o}rmander's vector fields and try to establish a Nash type result for weak solutions which will play a crucial role for corresponding nonlinear problems.

Let $\Omega\subset\mathbb{R}^n$ be a bounded domain, ${\mathcal{Q}}_{T}=\Omega\times(0,T)\subset{{\mathbb{R}}^n}\times{\mathbb{R}}, T>0$. Consider the following divergence parabolic equation
\begin{equation}\label{eq11}
{u_t}+{X_j^*}\left({a^{ij}(x,t)}{X_i}u \right)+{b_i}(x,t){X_i}u+c(x,t)u=f(x,t)-{X_i^*}{f^i}(x,t), ~ (x,t)\in{\mathcal{Q}}_{T}
\end{equation}
where ${X_i}=\sum\limits_{k = 1}^n {b_{ik}(x)}\frac{\partial}{\partial x_k} ({b_{ik}(x)}\in {C^\infty}(\Omega), i = 1,...,q, q \leq n)$ is the smooth vector field, $ {X_j^*} = -\sum\limits_{k = 1}^n \frac{\partial}{\partial x_k}\left({b_{jk}(x)}\cdot\right)$ is the adjoint of $X_j,j=1,...,q $. The summation coonventions in (\ref{eq11}) are omitted. Denote $$Lu={X_j^*}\left({a^{ij}(x,t)}{X_i}u \right)+{b_i}(x,t){X_i}u+c(x,t)u,$$ then (\ref{eq11}) is written by $${u_t}+Lu=f(x,t)-{X_i^*}{f^i}(x,t).$$

 Throughout this paper we make the following assumptions:

(C1) ${a^{ij}(x,t)}\in{L^\infty}({\mathcal{Q}}_{T})(i,j=1,2,...q)$ and there exists $\Lambda>0$ such that $$
{\Lambda ^{ - 1}}{\left| \xi  \right|^2} \leqslant {a^{ij}}(x,t){\xi _i}{\xi _j} \leqslant \Lambda {\left| \xi  \right|^2},{\kern 1pt} {\kern 1pt} {\kern 1pt} {\kern 1pt} (x,t) \in {\kern 1pt} {\mathcal{Q}}_{T} {\kern 1pt} ,{\kern 1pt} {\kern 1pt} {\kern 1pt} \xi  \in {\mathbb{R}^q};$$

(C2) for $m>(Q+2)/2,Q$ is the local homogeneous dimension relative to $\Omega,{b_i}(x,t)(i=1,2,...,q)$ and $ c(x,t)$ satisfy
\[\sum\limits_{i}{\|{b^{i}}^2\|}_{L^{m}({\mathcal{Q}}_{T})}+{\|c\|}_{L^{m}({\mathcal{Q}}_{T})}\leqslant\Lambda;\]

(C3) $f\in{L^{\frac{p(Q+2)}{Q+2+p}}({\mathcal{Q}}_{T})}, {f^i}\in{L^{p}({\mathcal{Q}}_{T})}, p>Q+2, i = 1,...,q $.

Let $$(Lu,\varphi)=\int_{\Omega}[(a^{ij}{X_i}u){X_j}\varphi+({b_i}{X_i}u+cu)\varphi]dx, ~for~ \varphi\in {W_{2;0}^{1,1}({\mathcal{Q}}_{T})};$$  we call that $u$ is a weak sub-solution (super-solution) to the equation (\ref{eq11}), if $u\in{V_2}({\mathcal{Q}}_{T})$ satisfies
\begin{equation}\label{eq12}
\int_{0}^{t}({u_t},\varphi)d\tau+\int_{0}^{t}(Lu,\varphi)d\tau\leqslant(\geqslant)
 \int_{0}^{t}[(f,\varphi)+(f^i),{X_i}\varphi]d\tau, t\in(0,T),
\end{equation}
for any $\varphi\in{W_{2;0}^{1,1}}({\mathcal{Q}}_{T}),\varphi\geqslant0.$ If u is not only a sub-solution but also a super-solution to (\ref{eq11}), we say that $u$ is a weak solution to (\ref{eq11}); at this time,
(\ref{eq12}) becomes an integral equality and the restriction $\varphi\geqslant0$ is eliminated. Spaces ${W_{2;0}^{1,1}}({\mathcal{Q}}_{T}),V_{2}({\mathcal{Q}}_{T})$ here and $V_{2}^{1,0}({\mathcal{Q}}_{T})$ appeared in the following paragraph will be described in Section 2 in detail.

We call $u\in{V_{2}^{1,0}({\mathcal{Q}}_{T})}$ belongs to the De Giorgi class
$$DG({\mathcal{Q}}_{T}):=DG\left({\mathcal{Q}}_{T};\lambda_0,\eta,M,{F_0},\gamma(\cdot),\delta\right),$$
if ${\|u\|}_{L^{\infty}({\mathcal{Q}}_{T})}\leqslant M $ and for $0<ess\sup\limits_{\mathcal{Q}_{\rho,\tau}}(u-k)_{\pm}\leqslant{\delta M},\delta\in(0,1],$
\begin{align}\label{eq13}
&\max\left\{\sup\limits_{0<{t_0}<t<{t_0}+\tau}\|\zeta(u-k)_{\pm}(\cdot,t)\|_{L^2(B_\rho)}^2,
\lambda_0\|X(\zeta(u-k)_{\pm})\|_{L^2(\mathcal{Q}_{\rho,\tau})}^2\right\}\nonumber \\
&\leqslant(1+\epsilon)\|\zeta(u-k)_{\pm}(\cdot,t_0)\|_{L^2(B_\rho)}^2\nonumber \\
& ~~+{\gamma(\epsilon)}\left\{[\|X\zeta\|_{L^{\infty}(\mathcal{Q}_{\rho,\tau})}^{2}
+\|\zeta_t\|_{L^{\infty}(\mathcal{Q}_{\rho,\tau})}+\|\zeta\|_{L^{\infty}(\mathcal{Q}_{\rho,\tau})}^{2}]
\|(u-k)_{\pm}\|_{L^2(\mathcal{Q}_{\rho,\tau})}^{2}\right.\nonumber \\
& ~~~~~~~~~~ \left.+({k^2}+{F_0}^2)|\mathcal{Q}_{\rho,\tau}\cap[(u-k)_{\pm}>0]|^{1-2/\eta}\right\},
\end{align}
where ${{\mathcal{Q}}_{\rho,\tau}}={B_{\rho}}\times{({t_0},{t_0}+\tau]}\subset{\mathcal{Q}_T},0<\rho,\tau<1,
\zeta\in{W_{2;0}^{1,1}}(\mathcal{Q}_{\rho,\tau}),0\leqslant{\zeta(x,t)}\leqslant1,\epsilon\in(0,1],\eta>Q+2,{F_0}\geqslant0,
\lambda_0>0$ is parameter, $\gamma(\cdot)$ is a non-negative decreasing functions. If (\ref{eq13}) is holds for $(u-k)_{+}$ ((\ref{eq13}) is holds for $(u-k)_-$), we denote $u\in{DG}^{+}(\mathcal{Q}_T)\left( u\in {DG}^{-}{\mathcal{Q}_T}\right)$. Clearly,
\[{DG(\mathcal{Q}_T)}={DG}^{+}(\mathcal{Q}_T)\bigcap{DG}^{-}(\mathcal{Q}_T).\]

Now we state the main results of this paper.

\begin{theorem}[H\"{o}lder Regularity]\label{Th11}
Let $u\in{V_{2}^{1,0}({\mathcal{Q}}_{T})}$ with ${\|u\|}_{L^{\infty}({\mathcal{Q}}_{T})}\leqslant M $ be the weak solution to (\ref{eq11}) with $(C1),(C2)$ and $(C3)$, then for ${\mathcal{Q}}\subset\subset{{\mathcal{Q}}_{T}},$ there exist $C=C(Q,\eta, \Lambda,\delta,{d_{\mathcal{P}}})\geqslant1$, and $\beta,0<\beta\leqslant{1-\frac{Q+2}{\eta}}$, such that
\begin{equation}\label{eq14}
{[u]}_{\beta;{\mathcal{Q}}} \leqslant C{{\overline{d}}^{-\beta}\left(M+{F_0}{\overline{d}}{|{\mathcal{Q}}_{T}|}^{-\frac{1}{\eta}}\right)},
\end{equation}
where ${\overline{d}}=\min\{1,d_{\mathcal{P}}({\mathcal{Q}},{\partial}_{\mathcal{P}}{{\mathcal{Q}}_{T}})\}, d_{\mathcal{P}}$ is the parabolic metric (see (\ref{eq26}) below ) and
${F_0}={\sum\limits_{i}\|f^i\|_{{L^p}({\mathcal{Q}_T})}}+{\|f\|_{L^{\frac{p(Q+2)}{Q+2+p}}({\mathcal{Q}}_{T})}}$.
\end{theorem}

\begin{theorem}[Harnack Inequality]\label{Th12}
Let $u\in{V_{2}^{1,0}({\mathcal{Q}}_{T})}$ with ${\|u\|}_{L^{\infty}({\mathcal{Q}}_{T})}\leqslant M $ be the weak solution to (\ref{eq11}) with $(C1),(C2)$ and $(C3)$,and $u\geqslant0$ on ${{\mathcal{Q}}_{R}^{a'}}={B_{8R}(x_0)}\times{({t_0},{t_0}+{a'}{R^2}]}\subset{{\mathcal{Q}}_{T}},{a'}>1,$ then
\begin{equation}\label{eq15}
\inf\limits_{{B_{R/2}}u({x_0},{t_0}+{a'}{R^2})}\geqslant{C^{-1}u({x_0},{t_0}+R^2)-
C{F_0}{R}{|{\mathcal{Q}}_{R}|}^{-\frac{1}{\eta}}},
\end{equation}
where $C$  depends on $Q,\Lambda,\eta,$ and $({a'}-1)^{-1}$.
\end{theorem}

The proofs of Theorems are based on the readjustment of De Giorgi's approach, and some new ingredients applying to our setting are replenished. It is worth emphasizing that we do not impose any artifical condition to the measure of metric ball induced by H\"{o}rmander's vector fields.

We note that (\ref{eq11}) involves the special case $${u_t}+{X_j^*}\left({a^{ij}(x,t)}{X_i}u \right)=0.$$

The rest of the paper is organized as follows: Section 2 explains H\"{o}rmander's vector fields, the parabolic Campanato space and parabolic H\"{o}lder space; several new preliminary results including the parabolic Sobolev inequality, (1,1) type Poincar\'{e} inequality of H\"{o}rmander's vector fields and a De Giorgi type Lemma are inferred. In section 3 we prove that weak solutions to (\ref{eq11}) are actually in the De Giorgi class $ DG({\mathcal{Q}}_{T})$ and derive some properties of functions in $ DG({\mathcal{Q}}_{T})$. Section 4 is devoted to proofs of main results. Theorem \ref{Th11} is proved based on the estimates of oscillations of solutions and the isomorphism between the parabolic Campanato space and parabolic H\"{o}lder space. Theorem \ref{Th12} is followed by using Theorem \ref{Th11} and an extension property of positivity for functions in the De Giorgi class.

%%---------------------------------------------------------------------------------
\section{Preliminaries}\label{Section 2}
%%---------------------------------------------------------------------------------

Let ${X_1},...,{X_q}{\kern 1pt} {\kern 1pt}(q \le n)$ are ${C^\infty }$ vector fields in $\Omega\subset{\mathbb{R}^N}.$ Throughout this paper, we always suppose that these vector fields satisfy the finite rank condition \cite{H1}, i.e., there exists a positive integer $s$ such that $\left\{{X_{\beta}(x_0)}\right\}_{|\beta|\leqslant s}$ spans ${\mathbb{R}^N}$ at every point ${x_0}\in\Omega\subset{\mathbb{R}^N}.$
 % $$ rank Lie\left[{X_1},...,{X_q}\right]\equiv n. $$ %

\begin{definition}[Carnot-Carath\'{e}odory distance,\cite{FP}]\label{De21}
An absolutely continuous curve $\gamma :[0,T] \to \Omega $ is said sub-unitary, if it is Lipschitz continuous and satisfies that for every $\xi  \in {\mathbb{R}^N}$ and ${\rm{a}}{\rm{.e}}{\rm{. }} ~t \in [0,T],$
\[ < \gamma '(t),\xi{ > ^2}\leqslant \sum\limits_{j = 1}^q { < {X_j}(\gamma (t)),\xi{ > ^2}}.{\kern 1pt}{\kern 1pt}\]
Let $\Phi(x,y)$ be the class of sub-unitary curves connected x and y, we define the Carnot-Carath\'{e}odory distance
(C-C distance) by \[ d(x,y) = \inf \{ T\ge 0:\gamma  \in \Phi (x,y)\}.\]
\end{definition}

The C-C metric ball is defined by
\[{B_R}({x_0}) = B({x_0},R) = \{ x \in \Omega :d({x_0},x) < R\} \]
and the Lebesgue measure of metric ball ${B_R}(x_0)$  by $|{B_R}(x_0)|$ . A fundamental doubling property with respect to the metric balls was showed in \cite{NSW}, namely, there are positive constants $C_1>1$ and $R_0$, such that for ${x_0}\in\Omega$ and $0<R<{R_0}$,
\begin{equation}\label{eq21}
 |B({x_0},2R)|\leqslant{C_1}|B({x_0},R)|,
\end{equation}
where $Q=\log_{2}C_1$, $Q$ acts as a dimension and is called the local homogeneous dimension relative to $\Omega$. It is easy to see from (\ref{eq21}) that for any $0<R\leqslant {R_0}$ and $\theta\in(0,1)$,  .
\begin{equation}\label{eq22}
|B_{\theta R}|\geqslant {{C_1}^{-1}}\theta^Q|B_R|,
\end{equation}
where $C_1$ and $R_0$ are constants in (\ref{eq21}).

The gradient of $u\in{C^1}(\Omega)$ is denoted by $Xu=\left({X_1}u,...,{X_q}u\right)$, and the norm of $Xu$ is of the form $|Xu|=\left(\sum\limits_{j=1}^q({X_j}u)^2\right)^{1/2}$. The Sobolev space $S_0^{1,p}(\Omega)(1\leqslant p<\infty)$  related to vector fields ${X_1},...,{X_q}{\kern 1pt}$ is the completion of $C_0^1(\Omega)$ under the norm
\begin{equation}\label{eq23}
 \|u\|_{S_0^{1,p}(\Omega)}:=\left[\int_{\Omega}\left({|u|^p}+|Xu|^p\right)\right]^{\frac{1}{p}}.
\end{equation}

\begin{definition}\label{De22}
The parabolic Sobolev space ${V_2}(\mathcal{Q}_T)$ on vector fields ${X_1},...,{X_q}$ is the set of all functions $u$ satisfying $Xu \in {L^2}(\mathcal{Q}_T)$ and $\sup\limits_{t\in(0,T)}\int_{\Omega_T}{|u|^2}dx<\infty$. The norm on ${V_2}(\mathcal{Q}_T)$ is
\begin{equation}\label{eq24}
\|u\|_{V_2}:=\left(\int_{\mathcal{Q}_T}|Xu|^2dxdt+\sup\limits_{t\in(0,T)}\int_{\Omega_T}{|u|^2}dx\right)^{\frac{1}{2}}<\infty.
 \end{equation}
\end{definition}

In the sequel, we also need spaces $V_2^{1,0}(\mathcal{Q}_T)$ and $W_2^{1,1}(\mathcal{Q}_T)$, where $V_2^{1,0}(\mathcal{Q}_T)$ is the set of  functions in ${V_2}(\mathcal{Q}_T)$ satisfying
\[\lim\limits_{h\rightarrow 0}\|u(\cdot,t+h)-u(\cdot,t)\|_{{L^2}(\Omega)}=0,t,t+h\in [0,T];\]
$W_2^{1,1}(\mathcal{Q}_T)$ contains functions satisfying $u\in{{L^2}(\mathcal(Q)_T)},Xu\in{{L^2}(\mathcal{Q}_T)}$ and ${u_t}\in{{L^2}(\mathcal{Q}_T)}$. The norm on $W_2^{1,1}(\mathcal{Q}_T)$ is
\begin{equation}\label{eq25}
\|u\|_{W_2^{1,1}}:=\left(\int_{\mathcal{Q}_T}|u|^2 dxdt +\int_{\mathcal{Q}_T}|Xu|^2 dxdt
  +\int_{\mathcal{Q}_T}{u_t^2} dxdt\right)^{\frac{1}{2}}.
\end{equation}
Obviously, we have $$W_2^{1,1}(\mathcal{Q}_T)\subset V_2^{1,0}(\mathcal{Q}_T)\subset {V_2}(\mathcal{Q}_T). $$
Moreover, ${\mathop{V}\limits^\circ}_2(\mathcal{Q}_T)$ and $W_{2;0}^{1,1}(\mathcal{Q}_T)$ are collections of functions in ${V_2}(\mathcal{Q}_T)$ and $W_2^{1,1}(\mathcal{Q}_T)$ satisfying $${u(\cdot,t)}|_{\partial\Omega}=0, {\rm{a}}{\rm{.e}}{\rm{. }}~t\in(0,T)$$ respectively, where $\partial\Omega$ is the boundary of $\Omega$ .

Let ${\mathcal{Q}}\subset\subset{\mathcal{Q}_T}$, we define the parabolic metric $d_{\mathcal{P}}$:
\begin{equation}\label{eq26}
{d_{\mathcal{P}}\left((x,t),(y,s)\right)}=\left(d(x,y)^2 +|t-s|\right)^{1/2}, ~ for~ (x,t),(y,s)\in{\mathcal{Q}}.
\end{equation}

\begin{definition}[H\"{o}lder space]\label{De23}
For $\alpha\in (0,1]$, let $C^{\alpha}(\mathcal{Q})$ be the set of functions $ u:\mathcal{Q}\to \mathbb{R} $ satisfying
\begin{eqnarray*} \left[u\right]_{\alpha;\mathcal{Q}}:=\sup\left\{\frac{|u(x,t)-u(y,s)|}{d_{\mathcal{P}}^{\alpha}\left((x,t),(y,s)\right)},
(x,t),(y,s)\in \mathcal{Q}, {(x,t)}\neq {(y,s)} \right\} <\infty;
\end{eqnarray*}
its norm is $$\|u\|_{\alpha;\mathcal{Q}}=\left[u\right]_{\alpha;\mathcal{Q}}+\|u\|_{L^{\infty}(\mathcal{Q})}.$$
\end{definition}

\begin{definition}[Campanato space]\label{De24}
For $1\leqslant p < +\infty$ and $ \lambda \geqslant 0,$ if $u\in L^p(\mathcal{Q})$ satisfies
$$\left[u\right]_{p,\lambda}:=\left\{\sup\limits_{Z\in{\mathcal{Q}},0<R\leqslant d}\left(\frac{R^{-\lambda}}{\hat{\mathcal{Q}_R}}\right)\int_{\hat{\mathcal{Q}_R}}|u(Z)-u_{\hat{\mathcal{Q}_R}}|^p dZ\right\}^{\frac{1}{p}}<\infty, $$
where $ d={\rm {diam}} \mathcal{Q}, Z=(x,t),\hat{\mathcal{Q}_R}=({B_R}(x)\times(t,t+R^2])\cap\mathcal{Q},u_{\hat{\mathcal{Q}_R}}=
\frac{1}{|\hat{\mathcal{Q}_R}|}\int_{\hat{\mathcal{Q}_R}}u(Y)dY,Y=(y,s),$ then we say that $u$ belongs to the Campanato space $\mathcal{L}^{p,\lambda}(\mathcal{Q})$ with the norm $$ \|u\|_{p,\lambda}=\left[u\right]_{p,\alpha}+\|u\|_{L^p}. $$
\end{definition}

Let us state two useful cut-off functions $\xi(x)$ and $\eta (t)$ (\cite{DK},\cite{GN}) which satisfy
$ 0\leqslant \xi \leqslant 1, \xi =1 $ on $ B_\rho (B_\rho \subset B_R \subset \Omega, \xi=0$ outside $B_R$) and $|X\xi|\leqslant \frac{c}{R-\rho},$
\begin{eqnarray*}
 \eta(t)=\left\{ \begin{array}{cl}
 \frac{t-({t_0-R^2})}{{R^2}-{\rho^2}}, & t\in (t_0 - R^2,t_0 - \rho ^2),\\
 ~1,& t\in [t_0 - \rho^2,t_0 + R^2].\\
 \end{array} \right.
\end{eqnarray*}

The following result is well known.

\begin{lemma} [Sobolev inequality,\cite{CDG1,CDG2,Lu1,Lu2}]\label{Le21}
For $1 \leqslant p <Q ,$ there exist $C>0$ and $R_0 > 0$ such that for any $x\in\Omega$ and $0<R\leqslant R_0,$ we have for any $u\in{S_0^{1,p}(B_R)},B_R=B_R(x)$,
\begin{equation}\label{eq27}
{\left(\frac{1}{|B_R|}\int_{B_R}|u|^{pk}dx\right)}^{\frac{1}{pk}}\leqslant CR{\left(\frac{1}{|B_R|}\int_{B_R}|Xu|^{p}dx\right)}^{\frac{1}{p}},
\end{equation}
where $1\leqslant k \leqslant Q/(Q-p)$. In particular, let $k=\frac{Q}{Q-p}$ , then
\begin{equation}\label{eq28}
{\left(\int_{B_R}|u|^{pQ/(Q-p)}dx\right)}^{(Q-p)/pQ}\leqslant CR|B_R|^{-\frac{1}{Q}}{\left(\int_{B_R}|Xu|^{p}dx\right)}^{\frac{1}{p}}.
\end{equation}
\end{lemma}

In the light of (\ref{eq28}) we can prove

\begin{theorem}[Parabolic Sobolev inequality]\label{th21}
For $u\in {\mathop{V}\limits^\circ}_2 (\mathcal{Q}_T),$ it follows $u\in {L^{2(Q+2)/Q}(\mathcal{Q}_T)}$ and
\begin{equation}\label{eq29}
\int_{\mathcal{Q}_R}|u|^{2(Q+2)/Q}dxdt\leqslant C{R^2}|B_R|^{-\frac{2}{Q}}\max\limits_{t\in(0,T)}\left(\int_{B_R}|u(x,t)|^2 dx\right)^{2/Q}\int_{\mathcal{Q}_R}|Xu|^2dxdt,
\end{equation}
where $\mathcal{Q}_R={B_R}(x_0)\times(t_0 - R^2,t_0]\subset\mathcal{Q}_T,B_R=B_R(x_0)$ and
\begin{equation}\label{eq210}
\left\|u\right\|_{L^{2(Q+2)/Q}(\mathcal{Q}_R)}\leqslant CR|\mathcal{Q}_R|^{-\frac{1}{Q+2}}\|u\|_{{V_2}(\mathcal{Q}_T)}.
\end{equation}
\end{theorem}
\textbf{Proof.} Using the H\"{o}lder inequality and (\ref{eq28}) with $p=1$, it sees
\begin{align*}
\int_{B_R}|u|^{2(Q+2)/Q}dx &= \int_{B_R}|u|^{2/Q}|u|^{2(1+Q)/Q}dx \\ &\leqslant  \left(\int_{B_R}|u|^{2}dx\right)^{1/Q}\left[\int_{B_R}\left(|u|^{2(1+Q)/Q}\right)^{Q/(Q-1)}dx\right]^{(Q-1)/Q} \\
 &\leqslant {CR|B_R|^{-\frac{1}{Q}}}\left(\int_{B_R}|u|^{2}dx\right)^{1/Q}\int_{B_R}\left|X
 \left(u^{2(1+Q)/Q}\right)\right| dx \\ &\leqslant CR|B_R|^{-\frac{1}{Q}}\left(\int_{B_R}|u|^{2}dx\right)^{1/Q}\int_{B_R}\left|u\right|^{(2+Q)/Q}|Xu|dx \\ &\leqslant CR|B_R|^{-\frac{1}{Q}}\left(\int_{B_R}|u|^{2}dx\right)^{1/Q}\left(\int_{B_R}|u|^{2(2+Q)/Q}dx\right)^{1/2}
 {\left(\int_{B_R}{|Xu|^2} dx\right)^{1/2}}
\end{align*}
 and then
\begin{align*}
\int_{B_R}|u|^{2(Q+2)/Q}dx &\leqslant {CR|B_R|^{-\frac{2}{Q}}}\left(\int_{B_R}|u|^{2}dx\right)^{2/Q}\int_{B_R}|Xu|^2 dx \\&\leqslant CR|B_R|^{-\frac{2}{Q}}\left(\max\limits_s\int_{B_R}|u(x,s)|^{2}dx\right)^{2/Q}\int_{B_R}|Xu|^2 dx.
\end{align*}
Integrating it with respect to $t$, (\ref{eq29}) is derived.

We arrive at by (\ref{eq29}) and the Young inequality that
\begin{eqnarray*}
\left\|u\right\|_{L^{2(Q+2)/Q}(\mathcal{Q}_R)} &\leqslant& CR^{\frac{Q}{Q+2}}|B_R|^{-\frac{1}{Q+2}}\sup\limits_{0<t<T}\left\|u(\cdot,t)\right\|_{{L^2}(B_R)}^{\frac{2}{Q+2}}
 \left\|Xu\right\|_{{L^2}(\mathcal{Q}_R)}^{\frac{Q}{Q+2}}\\ &\leqslant&
CR\left(R^2|B_R|\right)^{-\frac{1}{Q+2}}\left(\sup\limits_{0<t<T}\left\|u(\cdot,t)\right\|_{{L^2}(\Omega)}+
 \left\|Xu\right\|_{{L^2}(\mathcal{Q}_T)}\right)\\
 &\leqslant& CR|\mathcal{Q}_R|^{-\frac{1}{Q+2}}
 \left(\sup\limits_{0<t<T}\left\|u(\cdot,t)\right\|_{{L^2}(\Omega)}^2+
 \left\|Xu\right\|_{{L^2}(\mathcal{Q}_T)}^2\right)^{\frac{1}{2}} \\
 &=& CR|\mathcal{Q}_R|^{-\frac{1}{Q+2}}\|u\|_{{V_2}(\mathcal(Q)_T)}
\end{eqnarray*}
and (\ref{eq210}) is proved.

\begin{lemma}\label{Le22}
For $0<\theta <1$ and $\kappa>1$, there exists a positive constant $C_\kappa$ such that
\[(1-\theta)\leqslant C_\kappa (1-\theta)^{1-\frac{1}{\kappa}}.\]
\end{lemma}
\textbf{Proof.} Denote $$f(\theta)=\frac{1-\theta}{1-{\theta}^{1-\frac{1}{\kappa}}},$$ then
$$\lim_{\theta\rightarrow 0^+} f(\theta)=1, \lim_{\theta\rightarrow 1^-} f(\theta)=\frac{\kappa}{\kappa -1}.$$
Let
\begin{eqnarray*}
 F(\theta)=\left\{ \begin{array}{cl}
  1, & \theta =0,\\
 f(\theta),& \theta\in(0,1),\\
 \frac{\kappa}{\kappa -1},& \theta =1,\\
 \end{array} \right.
\end{eqnarray*}
we have that $F(\theta)$ is continuous and uniformly bounded on $ [0,1] $, which implies $f(\theta)\leqslant C_{\kappa}$ for some ${C_\kappa}>0.$

\begin{theorem}[(1,1) type Poincar¨¦ inequality]\label{th22}
Let $u\in{W^{1,1}(B_R)}$ and $$ E_0=\{x\in{B_R}|u(x)=0\}.$$ If $|E_0|>0,$ then
\begin{equation}\label{eq211}
\int_{B_R}|u|dx\leqslant \frac{CR|B_R|}{|E_0|}\int_{B_R}|Xu|dx,
\end{equation}
where $C>0$ relies only on $Q$.
\end{theorem}
\textbf{Proof.} Consider first $u\in Lips(B(x,R))$, then by the result of Jerison \cite{J},
\[\int_{B_R}|u(x)-u_{B_R}|dx\leqslant CR\int_{B_R}|Xu|dx,\]
where $u_{B_R}=\frac{1}{|B_R|}\int_{B_R}u dx.$ Since the above inequality has the self-improvement property (see \cite{HK}), i.e. there exists $\kappa>1$, such that $$\left(\frac{1}{|B_R|}\int_{B_R}|u(x)-u_{B_R}|^\kappa dx\right)^{\frac{1}{\kappa}}\leqslant C \frac{R}{|B_R|}\int_{B_R}|Xu|dx, $$ it follows
\[\left(\int_{B_R}|u(x)-u_{B_R}|^\kappa dx\right)^{\frac{1}{\kappa}}\leqslant CR |B_R|^{\frac{1}{\kappa}-1}\int_{B_R}|Xu|dx.\] Applying
\begin{align*}
\left|u_{B_R}\right| &= \frac{1}{|B_R|}\left|\int_{B_R}u dx\right| \\
&\leqslant \frac{1}{|B_R|}\int_{{B_R}\setminus E_0}|u| dx \\
& \leqslant \frac{1}{|B_R|}\left(\int_{B_R}|u|^\kappa dx\right)^{\frac{1}{\kappa}}\left|{B_R}\setminus E_0\right|^{1-\frac{1}{\kappa}},
\end{align*}
it arrives at
\begin{eqnarray*}
\left(\int_{B_R}|u|^\kappa dx\right)^{\frac{1}{\kappa}} &\leqslant&
\left(\frac{1}{|B_R|}\int_{B_R}|u(x)-u_{B_R}|^\kappa dx\right)^{\frac{1}{\kappa}}+|u_{B-R}||B_R|^{\frac{1}{\kappa}} \\ &\leqslant & CR|B_R|^{\frac{1}{\kappa}-1}\int_{B_R}|Xu|dx+\frac{{(B_R)\setminus E_0}^{1-\frac{1}{\kappa}}}{|B_R|^{1-\frac{1}{\kappa}}}\left(\int_{B_R}|u|^\kappa dx\right)^{\frac{1}{\kappa}}
\end{eqnarray*}
and then
$$\left(1-\frac{{(B_R)\setminus E_0}^{1-\frac{1}{\kappa}}}{|B_R|^{1-\frac{1}{\kappa}}}\right)
\left(\int_{B_R}|u|^\kappa dx\right)^{\frac{1}{\kappa}} \leqslant CR|B_R|^{\frac{1}{\kappa}-1}\int_{B_R}|Xu|dx. $$
Noting
$$ 1-\frac{{(B_R)\setminus E_0}}{|B_R|}
\leqslant C_{\kappa}\left(1-\frac{{(B_R)\setminus E_0}^{1-\frac{1}{\kappa}}}{|B_R|^{1-\frac{1}{\kappa}}}\right)$$
we have from (\ref{Le22}) that
\[ C_\kappa ^{-1}\left(1-\frac{{(B_R)\setminus E_0}}{|B_R|}\right)\left(\int_{B_R}|u|^\kappa dx\right)^{\frac{1}{\kappa}}
\leqslant CR|B_R|^{\frac{1}{\kappa}-1}\int_{B_R}|Xu|dx. \]
Using it and the H\"{o}lder inequality, it follows
\begin{eqnarray*}
\frac{|E_0|}{|B_R|}\int_{B_R}|u|dx  & \leqslant &\frac{|E_0|}{|B_R|} |B_R|^{1-\frac{1}{\kappa}}
\left(\int_{B_R}|u|^\kappa dx\right)^{\frac{1}{\kappa}} \\ &\leqslant& CR|B_R|^{\frac{1}{\kappa}-1}|B_R|^{1-\frac{1}{\kappa}}\int_{B_R}|Xu|dx \\ &=& CR\int_{B_R}|Xu|dx.
\end{eqnarray*}
Now (\ref{eq211}) is obtained by combining it and the density of $Lip(B_R)$ in $W^{1,1}(B_R)$.

\begin{theorem}[De Giorgi type lemma]\label{Th23}
Let $u\in {W^{1,1}(B_R)}$ and $$ A(k)=\{x\in{B_R}|u(x)>k\}, for ~l>k, $$ we have
\begin{equation}\label{eq212}
\left(l-k\right)|A(l)| \leqslant {\frac{CR|B_R|}{|{B_R} \setminus A(k)|}}\int_{A(k)\setminus A(l)}|Xu|dx,
\end{equation}
where $C>0$ relies only on $Q$.
\end{theorem}
\textbf{Proof.} Denoting the function
\begin{eqnarray*}
 \hat{u}(x)=\left\{ \begin{array}{cl}
  ~l-k, & x\in A(l),\\
 u(x)-k,& x\in A(k)-A(l),\\
 ~~0,& x\in {B_R}-A(k),\\
 \end{array} \right.
\end{eqnarray*}
and noting $\hat{u}\in{W^{1,1}(B_R)}$ , it obtains by using (\ref{eq211}) that
\begin{eqnarray*}
\left(l-k\right)|A(l)|\leqslant \int_{B_R}|\hat{u}(x)|dx
\leqslant {\frac{CR|B_R|}{|{B_R}\setminus A(k)|}}\int_{{B_R}\setminus A(k)}|X\hat{u}|dx
 = C {\frac{CR|B_R|}{|{B_R}\setminus A(k)|}}\int_{A(k)\setminus A(l)}|Xu|dx.
 \end{eqnarray*}
 This proves (\ref{eq212}).

%%---------------------------------------------------------------------------------
\section{Several auxiliary lemmas \label{Section 3}}
%---------------------------------------------------------------------------------

For the weak sub-solution (super-solution) to (\ref{eq11}), we have

\begin{lemma}\label{Le31}
Let $u\in{V_{2}^{1,0}(\mathcal{Q}_T)}$ be the bounded weak sub-solution (or super-solution) to (\ref{eq11}) with $(C1),(C2)$ and $(C3)$ then
$$ u\in{DG^+}~(or ~ u\in{DG^-})$$  % \left({\mathcal{Q}}_{T};\lambda_0,\eta,M,{F_0},\gamma(\cdot),\delta\right)
where $p>Q+2, \eta=\min\left\{p,2m\right\}, \gamma(\cdot)$ relies only on $Q,\Lambda$ and $p$,
\begin{equation}\label{eq31}
F_0={\sum\limits_i}\|f^i\|_{{L^p}(\mathcal{Q}_T)} + \|f\|_{L^{\frac{p(Q+2)}{Q+2+p}(\mathcal{Q}_T)}}<\infty.
\end{equation}
\end{lemma}
\textbf{Proof.} We only prove the conclusion for the weak sub-solution and the proof for the weak super-solution is similar. Multiplying the test function ${\zeta^2}(u-k)_+$ to (\ref{eq11}) and integrating on ${B_\rho}\times(t_0,t)$, it yields
\begin{eqnarray*}
{\int_{t_0}^t}\left(u_t,\zeta^2(u-k)_+\right)dt &=& {\int_{t_0}^t}\int_{B_\rho}{u_t}\zeta^2{(u-k)_+} dx dt \\ &=& \frac{1}{2}\int_{B_\rho}\zeta^2{(u-k)_+}(\cdot,t)dx \mid_{t_0}^t -\frac{1}{2}\cdot2{\int_{t_0}^t}\int_{B_\rho}\zeta{\zeta}_t {(u-k)_+^2} dx dt \\ &=& \frac{1}{2}\int_{B_\rho}\zeta^2{(u-k)_+}(x,t)dx - \frac{1}{2}\int_{B_\rho}\zeta^2{(u-k)_+}(x,t_0)dx \\ & &~~~ -{\int_{t_0}^t}\int_{B_\rho}\zeta{\zeta}_t {(u-k)_+^2} dx dt
\end{eqnarray*} and so
\begin{align}\label{eq32}
& \frac{1}{2}\int_{B_\rho}\zeta^2{(u-k)_+}(x,t)dx - \frac{1}{2}\int_{B_\rho}\zeta^2{(u-k)_+}(x,t_0)dx \nonumber \\ & ~-{\int_{t_0}^t}\int_{B_\rho}\zeta{\zeta}_t {(u-k)_+^2} dxdt + {\int_{t_0}^t}\int_{B_\rho} a^{ij}{X_i}u{X_j}(\zeta^2(u-k)_+)dxdt \nonumber \\ & ~+ {\int_{t_0}^t}\int_{B_\rho}({b^i}{X_i}u+cu)(\zeta^2(u-k)_+)dxdt \nonumber \\ &\leqslant {\int_{t_0}^t}\int_{B_\rho} [{f^i}{X_i}({\zeta^2}(u-k)_+)+f{\zeta^2}(u-k)_+]dxdt,
\end{align}
where ${t_0}<t\leqslant {t_0}+\tau ,0<\tau<1$. Noticing ${X_i}u={X_i}(u-k)_+$ and
\begin{align*}
& {X_i}u{X_j}({\zeta^2}(u-k)_+)={X_i}u{X_j}({\zeta}\cdot\zeta(u-k)_+) \\ &= \zeta {X_i}u[(u-k)_+{X_j}\zeta+{X_j}(\zeta(u-k)_+)] \\ &= [{X_i}(\zeta(u-k)_+)-{(u-k)_+}{X_i}\zeta][(u-k)_+{X_j}\zeta+{X_j}(\zeta(u-k)_+)] \\ &= {X_i}(\zeta(u-k)_+){X_j}(\zeta(u-k)_+)-{(u-k)_+^2}{X_i}\zeta{X_j}\zeta,
\end{align*}
we obtain from (\ref{eq32}) that
\begin{align}\label{eq33}
&~ \frac{1}{2}\int_{B_\rho}\zeta^2{(u-k)_+}(x,t)dx + {\int_{t_0}^t}\int_{B_\rho}a^{ij}{X_i}u{X_j}({\zeta^2}(u-k)_+)dxdt \nonumber \\ &\leqslant \frac{1}{2}\int_{B_\rho}\zeta^2{(u-k)_+}(x,t_0)dx + {\int_{t_0}^t}\int_{B_\rho}\zeta{\zeta}_t {(u-k)_+^2} dxdt \nonumber \\ &~+{\int_{t_0}^t}\int_{B_\rho}{a^{ij}}{(u-k)_+^2}{X_i}\zeta{X_j}\zeta dxdt -{\int_{t_0}^t}\int_{B_\rho}\left({b^i}{X_i}u+cu\right){\zeta^2}{(u-k)_+} dxdt \nonumber \\ & ~~~~+{\int_{t_0}^t}\int_{B_\rho}[{f^i}{X_i}({\zeta^2}(u-k)_+)+f{\zeta^2}(u-k)_+]dxdt.
\end{align}
Denote the third, fourth and fifth term in the right hand side of (\ref{eq33}) respectively by
\begin{eqnarray*} \begin{gathered} I_1:={\int_{t_0}^t}\int_{B_\rho}{a^{ij}}{(u-k)_+^2}{X_i}\zeta{X_j}\zeta dx dt, \hfill \\  I_2:=-{\int_{t_0}^t}\int_{B_\rho} \left({b^i}{X_i}u+cu\right) {\zeta^2}{(u-k)_+} dx dt,   \hfill  \\  I_3:= {\int_{t_0}^t}\int_{B_\rho}[{f^i}{X_i}({\zeta^2}(u-k)_+)+f{\zeta^2}(u-k)_+] dx dt.  \\ \end{gathered}
\end{eqnarray*}
Applying (C1) to $I_1$, it has \[{I_1}\leqslant\Lambda{\int_{t_0}^t}\int_{B_\rho}{(u-k)_+^2}|X\zeta|^2 dxdt.\]
In virtue of the Cauchy inequality, $\left[u>k\right]=\left\{(x,t)\mid u>k \right\}$ and $$\zeta{X_i}u={X_i}(\zeta(u-k)_+)-{(u-k)_+}{X_i}\zeta,$$ we see
\begin{align*}
{I_2} &\leqslant {\int_{t_0}^t}\int_{B_\rho}\left|{b^i}{\zeta}{(u-k)_+}{\zeta}{X_i}u\right|dxdt + {\int_{t_0}^t}\int_{B_\rho}\left|\sqrt{|c|}\zeta{(u-k)_+}\sqrt{|c|}u\zeta\right|dxdt \\ &\leqslant {\gamma}(\epsilon)\int_{{\mathcal{Q}_{\rho,\tau}}\cap[u>k]}{\zeta}^2{(u-k)_+}^2\left[\sum (b^i)^2\right]dxdt +\epsilon \int_{{\mathcal{Q}_{\rho,\tau}}\cap[u>k]}[{X_i}(\zeta(u-k)_+)-{(u-k)_+}{X_i}\zeta]^2 dxdt \\ & ~~+\gamma(\epsilon)\int_{{\mathcal{Q}_{\rho,\tau}}\cap[u>k]}|c|{\zeta}^2{(u-k)_+}^2 dxdt + \epsilon \int_{{\mathcal{Q}_{\rho,\tau}}\cap[u>k]}|c|((u-k)+k)^2{\zeta^2}dxdt \\ &\leqslant {\gamma}(\epsilon)\int_{{\mathcal{Q}_{\rho,\tau}}\cap[u>k]}{\zeta}^2{(u-k)_+}^2\left[\sum (b^i)^2 +|c|\right]dxdt \\ & ~~+\epsilon \int_{{\mathcal{Q}_{\rho,\tau}}\cap[u>k]}[{X_i}(\zeta(u-k)_+)-{(u-k)_+}{X_i}\zeta]^2 dxdt+2 \epsilon \int_{{\mathcal{Q}_{\rho,\tau}}\cap[u>k]}|c|((u-k)^2+k^2){\zeta^2}dxdt \\ &\leqslant \gamma(\epsilon)\int_{{\mathcal{Q}_{\rho,\tau}}\cap[u>k]}\left\{{\zeta}^2{(u-k)_+}^2 \left[\sum (b^i)^2 +|c|\right] + |c|k^2\zeta^2 \right\}dxdt \\ &~~+2 \epsilon \int_{{\mathcal{Q}_{\rho,\tau}}\cap[u>k]}{\left|{X_i}(\zeta(u-k)_+)\right|}^2+{(u-k)_+^2}\left|X\zeta\right|^2 dxdt
\end{align*} and
\begin{align*}
{I_3} &\leqslant \gamma(\epsilon)\int_{{\mathcal{Q}_{\rho,\tau}}\cap[u>k]}\left\{\sum{(f^i)^2}+|f|\zeta^2(u-k)_+\right\}dxdt \\ &~~~ + \epsilon \int_{{\mathcal{Q}_{\rho,\tau}}\cap[u>k]}[{(u-k)_+}{X_j}\zeta+{X_j}(\zeta(u-k)_+)]^2 dxdt \\ &\leqslant \gamma(\epsilon)\int_{{\mathcal{Q}_{\rho,\tau}}\cap[u>k]}\left\{\sum{(f^i)^2}+|f|\zeta^2(u-k)_+\right\}dxdt \\ &~~~ + 2\epsilon \int_{{\mathcal{Q}_{\rho,\tau}}\cap[u>k]}{(u-k)_+^2}\left|X\zeta\right|^2+{\left|{X_i}(\zeta(u-k)_+)\right|}^2 dxdt.
\end{align*}
Substituting these estimates into (\ref{eq33}) and using $(C1)$ to the second term in the left hand side of (\ref{eq33}), it is not difficult to derive
\begin{align*}
&~ \frac{1}{2}\int_{B_\rho}\zeta^2{(u-k)_+}(x,t)dx + {\Lambda}^{-1}{\int_{t_0}^t}\int_{B_\rho}{\left|X(\zeta(u-k)_+)\right|}^2 dxdt \\ &\leqslant \frac{1}{2}\int_{B_\rho}\zeta^2{(u-k)_+}(x,t_0)dx + {\int_{t_0}^t}\int_{B_\rho}|{\zeta}_t| {(u-k)_+^2} dxdt +\Lambda {\int_{t_0}^t}\int_{B_\rho}{(u-k)_+^2}\left|X\zeta\right|^2 dxdt \\ &~~+\gamma(\epsilon)\int_{{\mathcal{Q}_{\rho,\tau}}\cap[u>k]}\left\{{\zeta}^2{(u-k)_+}^2 \left[\sum (b^i)^2 +|c|\right] + |c|k^2\zeta^2 \right\}dxdt \\ &~~~+\gamma(\epsilon)\int_{{\mathcal{Q}_{\rho,\tau}}\cap[u>k]}\left\{\sum{(f^i)^2}+|f|\zeta^2(u-k)_+\right\}dxdt \\ &~~~~+4\epsilon \int_{{\mathcal{Q}_{\rho,\tau}}\cap[u>k]}{(u-k)_+^2}\left|X\zeta\right|^2+{\left|{X_i}(\zeta(u-k)_+)\right|}^2 dxdt.
\end{align*}
Choosing $\epsilon=\frac{\Lambda^{-1}}{16}$ , we have
\begin{align}\label{eq34}
&~ \frac{1}{2}\int_{B_\rho}\zeta^2{(u-k)_+}(x,t)dx + \frac{3{\Lambda^{-1}}}{4}{\int_{t_0}^t}\int_{B_\rho}{|X(\zeta(u-k)_+)|}^2 dxdt \nonumber \\ &\leqslant \frac{1}{2}\int_{B_\rho}\zeta^2{(u-k)_+}(x,t_0)dx + C \left\{ {\int_{t_0}^t}\int_{B_\rho}{(u-k)_+^2}(|X\zeta|^2+|\zeta_t|)dxdt \right. \nonumber \\ &~~ +\int_{{\mathcal{Q}_{\rho,\tau}}\cap[u>k]}\left[{\zeta}^2{(u-k)_+^2} \left[{\sum {(b^i)^2}} +|c|\right] + |c|k^2\zeta^2\right]dxdt \nonumber \\ &~~~+\left. \int_{{\mathcal{Q}_{\rho,\tau}}\cap[u>k]}{\left[\sum{(f^i)^2}+|f|\zeta^2{(u-k)_+}\right]} dxdt \right\}.
\end{align}
Denote the third and fourth term in the right hand side of (\ref{eq34}) by
\begin{eqnarray*}\begin{gathered} {II_1}:= \int_{{\mathcal{Q}_{\rho,\tau}}\cap[u>k]}\left[{\zeta}^2{(u-k)_+^2} \left[{\sum {(b^i)^2}} +|c|\right] + |c|k^2\zeta^2\right]dx dt, \hfill \\ {II_2}:=\int_{{\mathcal{Q}_{\rho,\tau}}\cap[u>k]}{\left[\sum{(f^i)^2}+|f|\zeta^2{(u-k)_+}\right]} dx dt. \hfill
\end{gathered} \end{eqnarray*}
Employing (C2), (\ref{eq31}) and the H\"{o}lder inequality, it shows
\begin{align*}
{II_1} &\leqslant \left(\int_{{\mathcal{Q}_{\rho,\tau}}\cap[u>k]}\left[{\sum {(b^i)^2}} +|c|\right]^m dxdt\right)^{\frac{1}{m}}\left(\int_{{\mathcal{Q}_{\rho,\tau}}\cap[u>k]} \left[{\zeta}^2{(u-k)_+^2}\right]^{\frac{m}{m-1}} dxdt \right)^{1-\frac{1}{m}} \nonumber \\ &~~~~+k^2 \left(\int_{{\mathcal{Q}_{\rho,\tau}}\cap[u>k]}|c|^m dxdt\right)^{\frac{1}{m}}|{{\mathcal{Q}_{\rho,\tau}}\cap[u>k]}|^{1-\frac{1}{m}} \nonumber \\ &\leqslant \Lambda\left\|\zeta(u-k)_+\right\|_{2m/(m-1)}^2 + k^2 \Lambda |{{\mathcal{Q}_{\rho,\tau}}\cap[u>k]}|^{1-\frac{1}{m}},
\end{align*}
where $m>\frac{Q+2}{2}$. As $2<\frac{2m}{m-1}<\frac{2(Q+2)}{Q}$, it knows by the Interpolation inequality and the Young inequality that for $0<\theta<1$,
\begin{align*}
\left\|\zeta(u-k)_+\right\|_{2m/(m-1)}^2 &\leqslant \left\|\zeta(u-k)_+\right\|_{2(Q+2)/Q}^{2(1-\theta)} \left\|\zeta(u-k)_+\right\|_{2}^{2\theta} \\ &\leqslant \frac{\epsilon}{\Lambda}\left\|\zeta(u-k)_+\right\|_{2(Q+2)/Q}^{2} + \gamma(\epsilon) \left\|\zeta(u-k)_+\right\|_{2}^{2}.
\end{align*}
It implies from it and (\ref{eq210}) that
\begin{eqnarray*}
{II_1}\leqslant \epsilon \left\|\zeta(u-k)_+\right\|_{{V_2}(\mathcal{Q}_{\rho,\tau})}^{2} + \gamma(\epsilon) \left\|\zeta(u-k)_+\right\|_{2}^{2} + k^2 \Lambda |{{\mathcal{Q}_{\rho,\tau}}\cap[u>k]}|^{1-\frac{1}{m}}.
\end{eqnarray*}
Using (\ref{eq210}) to $II_2$, it gets
\begin{align*}
{II_2} &\leqslant \left(\int\limits_{{\mathcal{Q}_{\rho,\tau}}\cap[u>k]}{\left[\sum{(f^i)^2}\right]^{\frac{p}{2}}} dxdt\right)^{\frac{2}{p}}\left|{\mathcal{Q}_{\rho,\tau}}\cap[u>k]\right|^{1-\frac{2}{p}} \\ &~+\left(\int\limits_{{\mathcal{Q}_{\rho,\tau}}\cap[u>k]}|\zeta f|^{\frac{(Q+2)p}{Q+2+p}}dxdt \right)^{\frac{Q+2+p}{(Q+2)p}}\left(\int\limits_{{\mathcal{Q}_{\rho,\tau}}\cap[u>k]}|\zeta(u-k)_+|^{\frac{2(Q+2)}{Q}}dxdt \right)^{\frac{Q}{2(Q+2)}}\left|{\mathcal{Q}_{\rho,\tau}}\cap[u>k]\right|^{\frac{1}{2}-\frac{1}{p}}  \\ &\leqslant {F_0^2}\left|{\mathcal{Q}_{\rho,\tau}}\cap[u>k]\right|^{1-\frac{2}{p}} + \left\|\zeta f\right\|_{\frac{(Q+2)p}{Q+2+p}}\left\|\zeta(u-k)_+\right\|_{\frac{2(Q+2)}{Q}} \left|{\mathcal{Q}_{\rho,\tau}}\cap[u>k]\right|^{\frac{1}{2}-\frac{1}{p}} \\ &\leqslant {F_0^2}\left|{\mathcal{Q}_{\rho,\tau}}\cap[u>k]\right|^{1-\frac{2}{p}} + \epsilon \left\|\zeta(u-k)_+\right\|_{{V_2}(\mathcal{Q}_{\rho,\tau})}^2 + \gamma(\epsilon){F_0^2}\left|{\mathcal{Q}_{\rho,\tau}}\cap[u>k]\right|^{1-\frac{2}{p}}.
\end{align*}
Putting estimates for $II_1$ and $II_2$ into (\ref{eq34}) and choosing $\eta=\min\{2m,p\}$ ($m>\frac{Q+2}{2},p>Q+2$), it leads to
\begin{align*}
&~ \frac{1}{2}\int_{B_\rho}\zeta^2{(u-k)_+}(x,t)dx + \frac{3{\Lambda^{-1}}}{4}{\int_{t_0}^t}\int_{B_\rho}{|X(\zeta(u-k)_+)|}^2 dxdt \\ & \leqslant \frac{1}{2}\int_{B_\rho}\zeta^2{(u-k)_+}(x,t_0)dx + C \left\{ {\int_{t_0}^t}\int_{B_\rho}{(u-k)_+^2}(|X\zeta|^2+|\zeta_t|)dxdt \right. \\ &~~~+ 2\epsilon \left\|\zeta(u-k)_+\right\|_{{V_2}(\mathcal{Q}_{\rho,\tau})}^{2} + \gamma(\epsilon) \left\|\zeta(u-k)_+\right\|_{2}^{2} + k^2 \Lambda |{{\mathcal{Q}_{\rho,\tau}}\cap[u>k]}|^{1-\frac{1}{m}} \\ &~~~~\left. +{F_0^2}\left|{\mathcal{Q}_{\rho,\tau}}\cap[u>k]\right|^{1-\frac{2}{p}} + \gamma(\epsilon){F_0^2}\left|{\mathcal{Q}_{\rho,\tau}}\cap[u>k]\right|^{1-\frac{2}{p}}\right\}.
\end{align*}
Taking $\epsilon$ small enough and $\lambda_0 =\Lambda^{-1}$, we conclude
\begin{eqnarray*}
& &\sup\limits_{0<{t_0}<t<{t_0}+\tau}\|\zeta(u-k)_{\pm}(\cdot,t)\|_{L^2(B_\rho)}^2+
\lambda_0\|X(\zeta(u-k)_{\pm})\|_{L^2(\mathcal{Q}_{\rho,\tau})}^2\} \\ &\leqslant &\frac{1}{2}\int_{B_\rho}\zeta^2{(u-k)_+}(x,t_0)dx  \\ & &~~~+\gamma(\epsilon)\left\{\left[\|X\zeta\|_{{L^{\infty}}(\mathcal{Q}_{\rho,\tau})}^2 + \|\zeta_t\|_{{L^{\infty}}(\mathcal{Q}_{\rho,\tau})} +\|\zeta\|_{{L^{\infty}}(\mathcal{Q}_{\rho,\tau})}^2 \right]\|(u-k)_{\pm}\|_{L^2(\mathcal{Q}_{\rho,\tau})}^{2}\right. \\ & & ~~~~~~~~~~ \left.+({k^2}+{F_0}^2)|\mathcal{Q}_{\rho,\tau}\cap[(u-k)_{\pm}>0]|^{1-2/\eta}\right\}
\end{eqnarray*}
and the proof is completed.

\textbf{Remark 3.1} \textit { Lemma \ref{Le31} indicates that the weak solution to (\ref{eq11}) belongs to the De Giorgi class.}

Next we give several useful properties for functions in the De Giorgi class.

\begin{lemma}\label{Le32} Let $u\in{DG^+(\mathcal{Q}_T)}$,$$\mathcal{Q}_{2R}^a = B_{2R}(x_0)\times (t_0,t_0+aR^2]\subset\mathcal{Q}_T, ~for~ 0<R\leqslant\frac{1}{2},0<a\leqslant1,$$ $\mu\geqslant\ess~sup\limits_{\mathcal{Q}_{2R}^a}u.$ If $0<H:=\mu-k\leqslant\delta M,u $ satisfies
\begin{equation}\label{eq35}
|{B_R}(x_0)\cap[u(\cdot,t)>k]|\leqslant(1-\sigma)|B_R|,~for~0<\sigma<1,~t\in(t_0,t_0 + aR^2]\subset(0,T),
\end{equation}
then for any positive integer $s$, either
\begin{equation}\label{eq36}
H\leqslant {2^s}(M+F_0)R|\mathcal{Q}_R|^{-\frac{1}{\eta}},
\end{equation} or
\begin{equation}\label{eq37}
\left|\mathcal{Q}_R^a \cap\left[u>\mu-\frac{H}{2^s}\right]\right|\leqslant \frac{C}{\sigma\sqrt{as}}|\mathcal{Q}_R^a|,
\end{equation}
where $C$ relies only on $Q,\lambda_0,\eta,\delta$ and $\gamma(\cdot), \mathcal{Q}_\rho^a=B_\rho(x_0)\times(t_0,t_0+a\rho^2]$, for $0<\rho<2R$.
\end{lemma}
\textbf{Proof.} Denote
\begin{eqnarray*} \begin{array}{l} A_R(k,t)=B_R(x_0)\times \cap [u(\cdot,t)>k], \\ A_R(k)={\mathcal{Q}_\rho^a}\cap[u>k], \\ k_l=\mu-\frac{H}{2^l},~l=0,1,..., \end{array}
\end{eqnarray*}
then $k_l$ is increasing,$A_R(k_l)$ is decreasing, and $|A_R(k,t)|\leqslant(1-\sigma)|B_R|$ by (\ref{eq35}). Applying Theorem \ref{Th23} and (\ref{eq35}), we have for $t_0\leqslant t\leqslant t_0 +aR^2 $,
\begin{align*}
(k_{l+1}-k_l)^2|A_R(k_{l+1},t)|^2 &\leqslant \frac{CR^2|B_R|^2}{|{B_R}\backslash A_R(k_l ,t)|^2}\left(\int_{A_R(k_{l},t)\setminus{A_R(k_{l+1},t)}}|Xu| dx\right)^2 \\ &\leqslant \frac{CR^2|B_R|^2}{|{B_R}-(1-\sigma){B_R}|^2}\left(\int_{A_R(k_{l},t)\setminus{A_R(k_{l+1},t)}}|Xu| dx\right)^2 \\ &\leqslant \frac{CR^2}{\sigma^2}|{A_R(k_{l},t)\setminus{A_R(k_{l+1},t)}}|\int_{A_R(k_{l},t)}|Xu|^2 dx.
\end{align*}
Integrating it in $t$ over $(t_0,t_0 +aR^2)$ and noting $k_{l+1}-k_l = \mu-\frac{H}{2^{l+1}}-(\mu-\frac{H}{2^{l}})=\frac{H}{2^{l+1}}$ and
\begin{align*}
& \int_{t_0}^{t_0 +aR^2}|{A_R(k_{l},t)\setminus{A_R(k_{l+1},t)}}|dt \\ &=\int_{t_0}^{t_0 +aR^2}|B_R(x_0)\cap[k_l <u(\cdot,t) <k_{l+1}]|dt \\ &=|\mathcal{Q}_R^a\cap[k_l <u <k_{l+1}]|=|{A_R(k_{l})\setminus{A_R(k_{l+1})}}|,
\end{align*}
 it gets
\begin{align}\label{eq38}
& \int_{t_0}^{t_0 +aR^2}|{A_R(k_{l+1},t)}|dt \nonumber \\ &\leqslant \frac{2^{l+1}CR}{\sigma H} \int_{t_0}^{t_0 +aR^2}|{A_R(k_{l},t)\setminus{A_R(k_{l+1},t)}}|^{\frac{1}{2}}dt \left(\int_{A_R(k_{l},t)}|Xu|^2 dx \right)^{\frac{1}{2}}dt \nonumber \\ &\leqslant \frac{CR 2^{l}}{\sigma H} \left(\int_{t_0}^{t_0 +aR^2}|{A_R(k_{l},t)\setminus{A_R(k_{l+1},t)}}|dt\right)^{\frac{1}{2}}dt \left(\int_{t_0}^{t_0 +aR^2}\int_{A_R(k_{l},t)}|Xu|^2 dx dt \right)^{\frac{1}{2}} \nonumber \\ &=\frac{CR 2^{l}}{\sigma H}|{A_R(k_{l})\setminus{A_R(k_{l+1})}}|^{\frac{1}{2}}\left(\int_{\mathcal{Q}_R^a}|X(u-k)_+|^2 dxdt\right)^{\frac{1}{2}}.
\end{align}
Letting $\xi(x)$ be the cut-off function between $B_R(x_0)$ and $B_{2R}(x_0)$ and using $$|X\xi|^2+|\xi_t|+|\xi|^2 \leqslant\frac{C}{R^2},~~|u-k_l|\leqslant\mu-k_l=\frac{H}{2^l},$$ and (\ref{eq22}), we obtain from (\ref{eq13}) ( $k$ is changed to $k_l$) that
\begin{align}\label{e39}
& \int_{\mathcal{Q}_{2R}^a}|X(\xi(u-k)_+)|^2 dxdt \nonumber \\ &\leqslant (1+\epsilon)\left\|{\xi(u-k_l)_+} (\cdot,t_0)\right\|_{L^2(B_2R)}^2 \nonumber \\ &~~~+ \gamma(\epsilon)\left\{\left[\|X\xi\|_{{L^{\infty}}(\mathcal{Q}_{2R}^a)}^2 + \|\xi_t\|_{{L^{\infty}}(\mathcal{Q}_{2R}^a)} +\|\xi\|_{{L^{\infty}}(\mathcal{Q}_{2R}^a)}^2 \right]\|(u-k_l)_{+}\|_{L^2(\mathcal{Q}_{2R}^a)}^{2}\right. \nonumber \\ & ~~~~~ \left.+({k^2}+{F_0}^2)|\mathcal{Q}_{2R}^a|^{1-2/\eta}\right\} \nonumber \\ &\leqslant C \left(\frac{H^2 |B_R|}{4^l} + \frac{H^2|\mathcal{Q}_{R}^a|}{4^l R^2} + (M+F_0)^2|\mathcal{Q}_{R}^a|^{1-2/\eta}\right).
\end{align}

If (\ref{eq36}) is invalid, i.e., there exists $l, 0\leqslant l\leqslant s-1$, such that $(M+F_0)\leqslant 2^{-l}HR^{-1}|\mathcal{Q}_R|^{\frac{1}{\eta}}$, then we note $|\mathcal{Q}_R^a|\leqslant|\mathcal{Q}_R|=R^2|B_R|$ and so $|\mathcal{Q}_{R}^a|^{1-2/\eta}\leqslant |\mathcal{Q}_{R}|^{1-2/\eta}$ to arrive at from the previous estimate (\ref{e39}) that
\begin{align*}
& ~\int_{\mathcal{Q}_{R}^a}|X{(u-k)_+}|^2 dxdt \\ &\leqslant \int_{\mathcal{Q}_{2R}^a}|X(\xi(u-k)_+)|^2 dxdt \\ &\leqslant C\left(\frac{H^2 |B_R|}{4^l} + \frac{H^2|\mathcal{Q}_{R}^a|}{4^l R^2} + \frac{H^2|\mathcal{Q}_{R}|}{{4^l}R^2}^{2/\eta}|\mathcal{Q}_{R}^a|^{1-2/\eta}\right) \\ &\leqslant C \frac{H^2 |B_R|}{4^l}.
\end{align*}
Taking it into (\ref{eq38}), it obtains
\begin{align*}
|A_R(k_{l+1})| &=|\mathcal{Q}_R^a\cap[u>k_{l+1}]|=\int_{t_0}^{t_0 +aR^2}|B_R \cap [u(\cdot,t)>k_{l+1}]|dt \\ &= \int_{t_0}^{t_0 +aR^2}|{A_R(k_{l+1},t)}|dt \leqslant \frac{CR |B_R|^{1/2}}{\sigma }|{A_R(k_{l})\setminus{A_R(k_{l+1})}}|^{\frac{1}{2}}.
\end{align*}
Squaring both sides and summing with respect to $l$ from $0$ to $s-1$ , we have by using $|A_R(k_s)|\leqslant|A_R(k_{l+1})|, 0\leqslant l\leqslant s-1, |A_R(k_0)|\leqslant|\mathcal{Q}_R^a|$ and $aR^2|B_R|=|\mathcal{Q}_R^a|$ that
\begin{align*}
& ~~s|A_R(k_s)|^2 \leqslant \sum\limits_{l=0}^{s-1}|A_R(k_{l+1})|^2 \\ &\leqslant\sum\limits_{l=0}^{s-1}\frac{CR^2|B_R|}{\sigma^2 }|{A_R(k_{l})\setminus{A_R(k_{l+1})}}| \\ &\leqslant \frac{CR^2|B_R|}{\sigma^2 }|A_R(k_0)| \leqslant \frac{CaR^2|B_R|}{a\sigma^2}|\mathcal{Q}_R^a| \\ &=\frac{C}{a\sigma^2}|\mathcal{Q}_R^a|^2,
\end{align*}
which implies (\ref{eq37}).

\begin{lemma}\label{Le33}Let $u\in{DG^+(\mathcal{Q}_T)}$,$$\hat{\mathcal{Q}}_{2R} = B_{2R}(x_0)\times (t_0,t_0+R^2]\subset\mathcal{Q}_T, ~for~ 0<R\leqslant\frac{1}{2},$$ $\mu\geqslant\ess~sup\limits_{\hat{\mathcal{Q}}_{2R}^a}u,\hat{\mathcal{Q}}_{2R}^a = B_{2R}(x_0)\times (t_0,t_0+aR^2],0<a\leqslant 1. $ For $0<H:=\mu-k\leqslant\delta M, 0<\sigma<1, $ if $u$ satisfies
\begin{equation}\label{eq39}
|{B_R}(x_0)\cap[u(\cdot,t_0)>k]|\leqslant(1-\sigma)|B_R|,
\end{equation} then there exists a positive integer $s_0 =s_0(\sigma)\geqslant 1$ relying on $Q,\lambda_0,\eta,\delta,\sigma$ and $\gamma(\cdot)$, such that either
\begin{equation}\label{eq310}
H\leqslant {2^{s_0}}(M+F_0)R|\mathcal{Q}_R|^{-\frac{1}{\eta}},
\end{equation} or
\begin{equation}\label{eq311}
\sup\limits_{t_0 <t<t_0 +R^2} \left|B_R(x_0)\cap\left[u(\cap,t)>\mu-\frac{H}{2^{s_0}}\right]\right| \leqslant \left[ 1-\sigma +\frac{1}{2}\min\{\sigma,1-\sigma\} \right]|B_R|.
\end{equation}
\end{lemma}
\textbf{Proof.} Suppose that $\xi(x)$ is a cut-off function between $B_{\beta R}(x_0) (0<\beta<1)$ and $B_R(x_0)$, such that $|X\xi|^2+|\xi_t|+|\xi|^2\leqslant\frac{C}{(1-\beta)^2R^2}$, and denote
\begin{eqnarray*} \begin{gathered} \mathcal{Q}_{R}^a = B_{R}(x_0)\times (t_0,t_0+aR^2], 0<a\leqslant 1, \\ A_R^a(k)={\mathcal{Q}_R^a}\cap[u>k].  \end{gathered}
\end{eqnarray*}
We observe $ u-k\leqslant \mu-k =H $ and apply (\ref{eq13}) on $\mathcal{Q}_R^a $ to see
\begin{align}\label{eq312}
& \sup\limits_{t_0<t\leqslant t_0+aR^2} \left\|\xi(u-k)_{+}(\cdot,t) \right\|_{L^2(B_R)}^2 \nonumber \\ &\leqslant (1+\epsilon)\left\|\xi(u-k)_{+}(\cdot,t_0) \right\|_{L^2(B_R)}^2 +\gamma(\epsilon)\left\{\frac{CH^2}{(1-\beta)^2R^2}|A_R^a(k)|+(M+F_0)^2|\frac{C}{(1-\beta)^2R^2}|^{1-2/\eta}\right\}.
\end{align}
On the other hand, for any integer $s_1\geqslant 1 $, \[ u-k>\mu-\frac{H}{2^{s_1}}-k = H-\frac{H}{2^{s_1}}=(1-2^{-s_1})H, ~ on ~ B_{\beta R}\bigcap[u(\cdot,t)>\mu-\frac{H}{2^{s_1}}], \] it follows
\begin{align}\label{eq313}
\left\|\xi(u-k)_{+}(\cdot,t) \right\|_{L^2(B_R)}^2 &\geqslant \int_{B_{\beta R}}|u(\cdot,t)-k|^2 dx  \nonumber \\ &\geqslant \int_{B_{\beta R}\bigcap[u(\cdot,t)>\mu-\frac{H}{2^{s_1}}]}|u(\cdot,t)-k|^2 dx  \nonumber \\ &\geqslant (1-2^{-s_1})^{2}{H^2}\left|{B_{\beta R}\bigcap[u(\cdot,t)>\mu-\frac{H}{2^{s_1}}]}\right|.
\end{align}
If (\ref{eq310}) is invalid, i.e. $(M+F_0)\leqslant 2^{-s_1}HR^{-1}|\mathcal{Q}_R|^{\frac{1}{\eta}},$ we note (\ref{eq39}) and $$ \left\|\xi(u-k)_{+}(\cdot,t_0) \right\|_{L^2(B_R)}^2\leqslant{H^2}\left|{B_{R}(x_0)\bigcap[u(\cdot,t_0)>k]}\right|, $$ and obtain from (\ref{eq312}) and (\ref{eq313}) that
\begin{align}\label{eq314}
 & \sup\limits_{t_0<t\leqslant t_0+aR^2}\left|{B_{\beta R}\bigcap[u(\cdot,t)>\mu-\frac{H}{2^{s_1}}]}\right| \nonumber \\ &\leqslant \frac{1}{(1-2^{-s_1})^2{H^2}}\sup\limits_{t_0<t\leqslant t_0+aR^2} \left\|\xi(u-k)_{+}(\cdot,t) \right\|_{L^2(B_R)}^2 \nonumber \\ &\leqslant \frac{1}{(1-2^{-s_1})^2{H^2}} \left\{(1+\epsilon)\left|{B_{R}(x_0)\bigcap[u(\cdot,t_0)>k]}\right|\right. \nonumber \\ & ~~~~~~ +\left.\gamma(\epsilon)\left[ \frac{CH^2}{(1-\beta)^2R^2}|A_R^a(k)|+(M+F_0)^2|A_R^a(k)|^{1-2/\eta} \right]\right\} \nonumber \\ &\leqslant \frac{(1+\epsilon)(1-\sigma)}{(1-2^{-s_1})^2}|B_R| + C\gamma(\epsilon)\left[\frac{|A_R^a(k)|}{(1-\beta)^2R^2} + \frac{|A_R^a(k)|^{1-2/\eta}}{R^2|\mathcal{Q}_R|^{-2/\eta}}\right] \nonumber \\ &\leqslant \frac{(1+\epsilon)(1-\sigma)}{(1-2^{-s_1})^2}|B_R| + C\gamma(\epsilon)|B_R| \left[ \frac{1}{(1-\beta)^2}\frac{|A_R^a(k)|}{|\mathcal{Q}_R|} + \left(\frac{|A_R^a(k)|}{|\mathcal{Q}_R|}\right)^{1-2/\eta}\right],
\end{align}
where the fact $R^2|B_R|=|\mathcal{Q}_R|$ is used. Obviously, $\frac{\epsilon(1-\sigma)}{(1-2^{-s_1})^2}\leqslant 4\epsilon $ from $1-\sigma \leqslant 1$ and $\mathord{(1-2^{-s_1})^2\geqslant\frac{1}{4}};$ $ \frac{1}{(1-\beta)^2}>1; |B_{\beta R}|\geqslant C_1^{-1}\beta^Q|B_R|$ and $|B_R|\backslash|B_{\beta R}|\leqslant(1-C_1^{-1}\beta^Q)|B_R|$ by (\ref{eq22}) and $\frac{|A_R^a(k)|}{|\mathcal{Q}_R|}\leqslant(\frac{|A_R^a(k)|}{|\mathcal{Q}_R|})^{1-\frac{2}{\eta}}$. Combining these and (\ref{eq314}), it derives that for $t_0<t\leqslant t_0+aR^2,$
\begin{align}\label{eq315}
& \left|B_R(x_0)\bigcap\left[u(\cdot,t)>\mu-\frac{H}{2^{s_1}}\right]\right| \nonumber \\ &\leqslant |B_R|\backslash|B_{\beta R}| + \left|B_{\beta R}(x_0)\bigcap\left[u(\cdot,t)>\mu-\frac{H}{2^{s_1}}\right]\right| \nonumber \\ &\leqslant (1-C_1^{-1}\beta^Q)|B_R| + \frac{(1+\epsilon)(1-\sigma)}{(1-2^{-s_1})^2}|B_R|+ C\gamma(\epsilon)|B_R| \left[ \frac{1}{(1-\beta)^2}\frac{|A_R^a(k)|}{|\mathcal{Q}_R|} + \left(\frac{|A_R^a(k)|}{|\mathcal{Q}_R|}\right)^{1-2/\eta}\right] \nonumber \\ &\leqslant |B_R|\left[1-C_1^{-1}\beta^Q +\frac{(1-\sigma)}{(1-2^{-s_1})^2}+4\epsilon +\frac{C\gamma(\epsilon)}{(1-\beta)^2}\left(\frac{|A_R^a(k)|}{|\mathcal{Q}_R|}\right)^{1-2/\eta} \right].
\end{align}
Choosing $\beta\in(0,1)$ such that $$1-\beta =\left(\frac{|A_R^a(k)|}{|\mathcal{Q}_R|}\right)^{(1/3)(1-2/\eta)}$$ and watching $1-C_1^{-1}\beta^Q\leqslant {C_Q}(1-\beta)$ for a positive constant $C_Q$, it obtains from (\ref{eq315}) that
\begin{align}\label{eq316}
& \left|B_R(x_0)\bigcap\left[u(\cdot,t)>\mu-\frac{H}{2^{s_1}}\right]\right| \nonumber \\ &\leqslant |B_R|\left[{C_Q}\left(\frac{|A_R^a(k)|}{|\mathcal{Q}_R|}\right)^{(1/3)(1-2/\eta)} +\frac{(1-\sigma)}{(1-2^{-s_1})^2}+4\epsilon +C\gamma(\epsilon)\left(\frac{|A_R^a(k)|}{|\mathcal{Q}_R|}\right)^{(1/3)(1-2/\eta)}\right] \nonumber \\ &\leqslant |B_R|\left[\frac{(1-\sigma)}{(1-2^{-s_1})^2}+4\epsilon +({C_Q} +C\gamma(\epsilon))\left(\frac{|A_R^a(k)|}{|\mathcal{Q}_R|}\right)^{(1/3)(1-2/\eta)}\right].
\end{align}
Since $\gamma(\epsilon)$ is decreasing, we know that $ q=\epsilon[\gamma(\epsilon)]^{-1},\epsilon\in\mathbb{R}^+ $, is strictly increasing and its inverse function $\epsilon=\varphi(q)$ satisfies $\varphi(q)\to 0$ as $q\to 0$. Using $\gamma(\epsilon)=\frac{\epsilon}{q}$ and choosing $$\epsilon=\varphi\left(\left(\frac{|A_R^a(k)|}{|\mathcal{Q}_R|}\right)^{(1/3)(1-2/\eta)}\right),$$ it follows $$ \gamma(\epsilon)=\epsilon\left(\frac{|A_R^a(k)|}{|\mathcal{Q}_R|}\right)^{-(1/3)(1-2/\eta)} = \varphi\left(\left(\frac{|A_R^a(k)|}{|\mathcal{Q}_R|}\right)^{(1/3)(1-2/\eta)}\right) \left(\frac{|A_R^a(k)|}{|\mathcal{Q}_R|}\right)^{-(1/3)(1-2/\eta)} $$ and so (\ref{eq316}) becomes
\begin{align}\label{eq317}
& \left|B_R(x_0)\cap\left[u(\cap,t)>\mu-\frac{H}{2^{s_1}}\right]\right| \nonumber \\ &\leqslant |B_R|\left\{\frac{(1-\sigma)}{(1-2^{-s_1})^2} + C\varphi\left(\left(\frac{|A_R^a(k)|}{|\mathcal{Q}_R|}\right)^{(1/3)(1-2/\eta)}\right)\right\}.
\end{align}
As $\epsilon=\varphi(q)$ is increasing and $|A_R^a(k)|\leqslant|\mathcal{Q}_R^a|=a|\mathcal{Q}_R|,$ it implies $$C\varphi\left(\frac{|A_R^a(k)|}{|\mathcal{Q}_R|}\right)^{(1/3)(1-2/\eta)}\leqslant C\varphi(a^{(1/3)(1-2/\eta)}).$$
Picking $a=a(\sigma)>0$ satisfying $$ C\varphi(a^{(1/3)(1-2/\eta)})\leqslant\frac{1}{4}\min\{1-\sigma,\sigma\} $$
and using $\lim\limits_{s_1\to\infty}\frac{(1-\sigma)}{(1-2^{-s_1})^2}=1-\sigma$, we can take $s_1=s_1(\sigma)$ large enough, such that $$\frac{(1-\sigma)}{(1-2^{-s_1})^2}\leqslant 1-\sigma +\frac{1}{4}\min\{1-\sigma,\sigma\}.$$
With it and (\ref{eq317}), it derives
\begin{equation}\label{eq318}
\sup\limits_{t_0\leqslant t\leqslant t_0+aR^2}\left|B_R(x_0)\bigcap\left[u(\cdot,t)>\mu-\frac{H}{2^{s_1}}\right]\right|\leqslant\left\{ 1-\sigma +\frac{1}{2}\min\{1-\sigma,\sigma\}\right\}|B_R|.
\end{equation}
Similarly, denoting
\begin{eqnarray*} \begin{array}{c} \mathcal{Q}_{R}^{1-a} = B_{R}(x_0)\times (t_0+aR^2,t_0+R^2],  \\ A_R^{1-a}(k)={\mathcal{Q}_R^{1-a}}\bigcap[u>k],  \end{array}
\end{eqnarray*}
and repeating the process above, we have for a large $s_2\geqslant1,$
\begin{equation}\label{eq319}
\sup\limits_{t_0+aR^2\leqslant t\leqslant t_0+R^2}\left|B_R(x_0)\bigcap\left[u(\cdot,t)>\mu-\frac{H}{2^{s_2}}\right]\right|\leqslant\left\{ 1-\sigma +\frac{1}{2}\min\{1-\sigma,\sigma\}\right\}|B_R|.
\end{equation}
Letting $s_0=\max\{s_1,s_2\}$, it shows by combining (\ref{eq318}) and (\ref{eq319}) that
\[ \sup\limits_{t_0\leqslant t\leqslant t_0+R^2}\left|B_R(x_0)\bigcap\left[u(\cdot,t)>\mu-\frac{H}{2^{s_2}}\right]\right|\leqslant\left\{ 1-\sigma +\frac{1}{2}\min\{1-\sigma,\sigma\}\right\}|B_R|, \]
which is (\ref{eq311}).

\begin{lemma}\label{Le34} Let $u\in{DG^+(\mathcal{Q}_T)},\mathcal{Q}_{R} = B_{R}(x_0)\times (t_0-R^2,t_0]\subset\mathcal{Q}_T, 0<R\leqslant 1 $, and $\mu\geqslant\ess~sup\limits_{\mathcal{Q}_{R}} u. $ there exists $\theta\in(0,1)$ depending on $Q,\lambda_0,\eta,\delta,\sigma$ and $\gamma(\cdot),$  such that for $k<\mu,$ if
\begin{equation}\label{eq320}
|\mathcal{Q}_R\cap[u>k]|\leqslant\theta|\mathcal{Q}_R|,
\end{equation}
\begin{equation}\label{eq321}
\delta M \geqslant H:=\mu-k > (M+F_0)R|\mathcal{Q}_R|^{-\frac{1}{\eta}},
\end{equation} then
\begin{equation}\label{eq322}
\ess~sup\limits_{\mathcal{Q}_{R/2}} u\leqslant \mu-\frac{H}{2}.
\end{equation}
\end{lemma}
To prove Lemma \ref{Le34}, we recall a known result.

\begin{lemma}[\cite{C}]\label{Le35}
Given a non-negative sequence $\{y_h\}~(h=0,1,2,...)$  satisfying the recursion relation $$y_{h+1}\leqslant C{b^h}y_h^{1+\epsilon},$$ where $b>1$ and $\epsilon >0$, if $${y_0}\leqslant {\theta}=C^{-1/\epsilon}b^{-1/{\epsilon^2}},$$ then \[\lim_{h\to \infty}{y_h}=0.\]
\end{lemma}
\textbf{Proof of Lemma \ref{Le34}.} Denote
\begin{equation}\label{eq323}
{R_{\bar{m}}}=\frac{R}{2}+\frac{R}{2^{\bar{m}+1}},{k_{\bar{m}}}=\mu-\frac{H}{2}-\frac{H}{2^{\bar{m}+1}},~and~ \mathcal{Q}_{\bar{m}}=\mathcal{Q}_{R_{\bar{m}}}~\bar{m}=0,1,2,...,
\end{equation}
then $R_m$ is decreasing, $\mathcal{Q}_{\bar{m}}$ is increasing and $$ R_{\bar{m}}^2-R_{\bar{m+1}}^2 = \left(\frac{R^2}{4}+\frac{R^2}{2^{\bar{m}+1}}+\frac{R^2}{2^{2\bar{m}+2}}\right) - \left(\frac{R^2}{4}+\frac{R^2}{2^{\bar{m}+2}}+\frac{R^2}{2^{2\bar{m}+4}}\right)\geqslant \frac{R^2}{2^{\bar{m}+2}}.$$ Take a cut-off function $\zeta_{\bar{m}}(x,t)$ between $\mathcal{Q}_{\bar{m}}$ and $\mathcal{Q}_{\bar{m}+1}$, then $$ |X\zeta_{\bar{m}}|^2\leqslant\left(\frac{C}{R_{\bar{m}}-R_{\bar{m+1}}}\right)^2=\frac{C2^{2\bar{m}+4}}{R^2}, |\zeta_{\bar{m}t}|\leqslant\frac{C}{ R_{\bar{m}}^2-R_{\bar{m+1}}^2} \leqslant\frac{C2^{2\bar{m}+2}}{R^2},|\zeta_{\bar{m}}|\leqslant 1, $$ and $$ |X\zeta_{\bar{m}}|^2+|\zeta_{\bar{m}t}|+|\zeta_{\bar{m}}|^2\leqslant C\frac{2^{4\bar{m}}}{R^2}. $$
Using
\begin{align*}
\left\|{\zeta_{\bar{m}}}(u-k_{\bar{m}})_{+}(\cdot,t_0) \right\|_{L^2(B_{R_{\bar{m}}})}^2 &= \frac{1}{R^2}\int_{t_0-R^2}^{t_0}{\int_{B_{R_{\bar{m}}}}}[(u-k_{\bar{m}})_{+}(\cdot,t_0)]^2 dxdt \\ &\leqslant \left\|(u-k_{\bar{m}})_{+}\right\|_{L^2(\mathcal{Q}_{\bar{m}})}^2,
\end{align*}
and replacing $\mathcal{Q}_{\rho,\tau},\zeta$ and $k$ in (\ref{eq13}) by $\mathcal{Q}_{\bar{m}},\zeta_{\bar{m}}$ and $k_{\bar{m}}$, we obtain
\begin{align}\label{eq324}
\left\|{\zeta_{\bar{m}}}(u-k_{\bar{m}})_{+}\right\|_{V_2(\mathcal{Q}_{\bar{m}})}^2 &= \sup\limits_{t_0<t<t_0+\tau}\left\|\zeta_{\bar{m}}(u-k_{\bar{m}})_{+}(\cdot,t) \right\|_{L^2(B_{R_{\bar{m}}})}^2 + \left\|X(\zeta_{\bar{m}}(u-k_{\bar{m}})_{+})\right\|_{L^2(\mathcal{Q}_{\bar{m}})}^2 \nonumber \\ &\leqslant C\left[\frac{2^{4\bar{m}}}{R^2}\left\|(u-k_{\bar{m}})_{+}\right\|_{L^2(\mathcal{Q}_{\bar{m}})}^2 + (M+F_0)^2 |\mathcal{Q}_{\bar{m}}\bigcap[u>k_{\bar{m}}]|^{1-\frac{1}{\eta}} \right].
\end{align}
Denoting $A_{\bar{m}}=\mathcal{Q}_{\bar{m}}\cap[u>k_{\bar{m}}]$ it follows $\left\|(u-k_{\bar{m}})_{+}\right\|_{L^2(\mathcal{Q}_{\bar{m}})}^2\leqslant H^2|A_{\bar{m}}|$ from $u-{k_{\bar{m}}}\leqslant\mu-{k_{\bar{m}}}\leqslant H$. Applying it into (\ref{eq324}), we have by (\ref{eq210}),
\begin{align}\label{eq325}
& \left\|{\zeta_{\bar{m}}}(u-k_{\bar{m}})_{+}\right\|_{\frac{2(Q+2)}{Q},(\mathcal{Q}_{\bar{m}})}^2 \nonumber \\ &\leqslant C{R^2}|\mathcal{Q}_R|^{-\frac{2}{Q+2}} \left\|{\zeta_{\bar{m}}}(u-k_{\bar{m}})_{+}\right\|_{V_2(\mathcal{Q}_{\bar{m}})}^2 \nonumber \\ &\leqslant C|\mathcal{Q}_{R}|^{-\frac{2}{Q+2}}\left[2^{4{\bar{m}}}{H^2}|A_{\bar{m}}| + R^2 (M+F_0)^2 |A_{\bar{m}}|^{1-2/\eta}\right].
\end{align}
On the other hand, $u>k_{\bar{m}+1}$ on $A_{\bar{m}+1}$ and
\begin{align*}
\left\|{\zeta_{\bar{m}}}(u-k_{\bar{m}})_{+}\right\|_{\frac{2(Q+2)}{Q},\mathcal{Q}_{\bar{m}}}^2 &= \left(\int_{\mathcal{Q}_{R}} \left({\zeta_{\bar{m}}}(u-k_{\bar{m}})_{+}\right)^{\frac{2(Q+2)}{Q}}dxdt \right)^{\frac{Q}{Q+2}} \\ &\geqslant \left(\int_{A_{\bar{m}+1}} \left((u-k_{\bar{m}})_{+}\right)^{\frac{2(Q+2)}{Q}}dxdt \right)^{\frac{Q}{Q+2}} \\ &\geqslant ({k_{\bar{m}+1}}-{k_{\bar{m}}})^2 |{A_{\bar{m}+1}}|^{\frac{Q}{Q+2}},
\end{align*} it shows
\begin{equation}\label{eq326}
({k_{\bar{m}+1}}-{k_{\bar{m}}})^2 |{A_{\bar{m}+1}}| \leqslant \left\|{\zeta_{\bar{m}}}(u-k_{\bar{m}})_{+}\right\|_{\frac{2(Q+2)}{Q},(\mathcal{Q}_{R})}^2 |{A_{\bar{m}+1}}|^{\frac{2}{Q+2}}.
\end{equation}
Substituting ${k_{\bar{m}+1}}-{k_{\bar{m}}}=\frac{H}{2^{\bar{m}+2}}$ into (\ref{eq326}), we derive from (\ref{eq325}) that
\begin{align}\label{eq327}
|{A_{\bar{m}+1}}| &\leqslant \frac{2^{2{\bar{m}}+4}}{H^2}C|\mathcal{Q}_{R}|^{-\frac{2}{Q+2}}\left[2^{4{\bar{m}}}{H^2}|A_{\bar{m}}| + R^2(M+F_0)^2 |A_{\bar{m}}|^{1-2/\eta}\right]|{A_{\bar{m}+1}}|^{\frac{2}{Q+2}} \nonumber \\ &\leqslant 2^{4\bar{m}}C|\mathcal{Q}_{R}|^{-\frac{2}{Q+2}}\left[2^{4{\bar{m}}}|A_{\bar{m}}| + \frac{R^2(M+F_0)^2}{H^2} |A_{\bar{m}}|^{1-2/\eta}\right]|{A_{\bar{m}+1}}|^{\frac{2}{Q+2}} \nonumber \\ &\leqslant C 2^{8\bar{m}}|\mathcal{Q}_{R}|^{-\frac{2}{Q+2}}\left[|A_{\bar{m}}| + \frac{R^2(M+F_0)^2}{H^2} |A_{\bar{m}}|^{1-2/\eta}\right]|{A_{\bar{m}+1}}|^{\frac{2}{Q+2}}.
\end{align}
Due to (\ref{eq321}) and $|A_{\bar{m}}|\leqslant|\mathcal{Q}_R|$, it yields from (\ref{eq327}) that
\begin{align*}
|{A_{\bar{m}+1}}| &\leqslant C 2^{8\bar{m}}|\mathcal{Q}_{R}|^{-\frac{2}{Q+2}}\left(|A_{\bar{m}}| + |A_{\bar{m}}|^{1-\frac{2}{\eta}}|\mathcal{Q}_R|^{\frac{2}{\eta}}\right)|{A_{\bar{m}+1}}|^{\frac{2}{Q+2}} \\ &=C 2^{8\bar{m}}|\mathcal{Q}_{R}|^{-\frac{2}{Q+2}}|A_{\bar{m}}|^{1-\frac{2}{\eta}}\left(|A_{\bar{m}}|^{\frac{2}{\eta}} + |\mathcal{Q}_R|^{\frac{2}{\eta}}\right)|{A_{\bar{m}+1}}|^{\frac{2}{Q+2}} \\ &\leqslant C2^{8\bar{m}}|A_{\bar{m}}|^{1-\frac{2}{\eta}+\frac{2}{Q+2}}|\mathcal{Q}_R|^{\frac{2}{\eta}-\frac{2}{Q+2}},
\end{align*} hence
\begin{equation}\label{eq328}
\frac{|A_{\bar{m}+1}|}{|\mathcal{Q}_R|}\leqslant C2^{8\bar{m}}\left(\frac{|A_{\bar{m}}|}{|\mathcal{Q}_R|}\right)^{1-\frac{2}{\eta}+\frac{2}{Q+2}}.
\end{equation}

Let $$y_{\bar{m}}=\frac{|A_{\bar{m}}|}{|\mathcal{Q}_R|},$$ then (\ref{eq328}) becomes $$y_{\bar{m}+1}\leqslant C2^{8\bar{m}}y_{\bar{m}}^{1-\frac{2}{\eta}+\frac{2}{Q+2}}. $$
Observing $1-\frac{2}{\eta}+\frac{2}{Q+2}>1$ (as $\eta>Q+2$ ) and (\ref{eq320}), we have
\begin{equation}\label{eq329}
y_0=\frac{|A_0|}{|\mathcal{Q}_R|}=\frac{|\mathcal{Q}_{R}\cap[u>k]|}{|\mathcal{Q}_R|}\leqslant \theta,~\theta\in(0,1),
\end{equation}
then $\lim\limits_{\bar{m}\to\infty}y_{\bar{m}}=0$ by Lemma \ref{Le35}, therefore $\lim\limits_{\bar{m}\to\infty}|A_{\bar{m}}|=0.$

Since $$ R_{\bar{m}}\to\frac{R}{2},k_{\bar{m}}\to\mu-\frac{H}{2} ~as~ m\to\infty $$ we have
\[ 0=\lim_{\bar{m}\to\infty}|A_{\bar{m}}|=\lim_{\bar{m}\to\infty}|{\mathcal{Q}_{R_{\bar{m}}}}\cap[u>k_{\bar{m}}]|= |\mathcal{Q}_{R/2}\cap[u>\mu-\frac{H}{2}]|, \] which gives (\ref{eq322}).

\begin{lemma}\label{Le36}  Let $u\in{DG^+(\mathcal{Q}_T)} (u\in{DG^-(\mathcal{Q}_T)})$,
$$\hat{\mathcal{Q}}_{2R} = B_{2R}(x_0)\times (t_0-R^2,t_0]\subset\hat{\mathcal{Q}}_T, ~for~ 0<R\leqslant 1/2,$$
$$ \mu\geqslant \ess~sup \limits_{\hat{\mathcal{Q}}_{2R}}u ~(\tilde{\mu}\leqslant ess\inf \limits_{\hat{\mathcal{Q}}_{2R}}u). $$
If for $0<\mu-k\leqslant\delta M ~(0<k-\tilde{\mu}\leqslant\delta M)$ and $0<\sigma<1,u $ satisfies
\begin{equation}\label{eq330}
|{B_R}(x_0)\cap[u(\cdot,t_0-R^2)\geqslant k]|\leqslant(1-\sigma)|B_R| ~ (|{B_R}(x_0)\cap[u(\cdot,t_0-R^2)\leqslant k]|\leqslant(1-\sigma)|B_R|),
\end{equation} then there exists $s=s(\sigma)\geqslant 1$ depending on $Q,\lambda_0,\eta,\delta,\sigma$ and $\gamma(\cdot)$ such that either
\begin{equation}\label{eq331}
H:=\mu-k\leqslant {2^{s}}(M+F_0)R|\mathcal{Q}_R|^{-\frac{1}{\eta}} ~(H:=k-\tilde{\mu}\leqslant {2^{s}}(M+F_0)R|\mathcal{Q}_R|^{-\frac{1}{\eta}}),
\end{equation} or
\begin{equation}\label{eq332}
\ess~sup\limits_{\mathcal{Q}_{R/2}} u\leqslant \mu-\frac{H}{2^s} ~(ess\inf\limits_{\mathcal{Q}_{R/2}} u\geqslant \tilde{\mu}+\frac{H}{2^s}),
\end{equation}
where $\mathcal{Q}_{\rho} = B_{\rho}(x_0)\times (t_0-\rho^2,t_0].$
\end{lemma}
\textbf{Proof.} Let $s_0$ be the constant in Lemma \ref{Le33}. If (\ref{eq331}) is invalid for $s\geqslant s_0$, then from Lemma \ref{Le33} and (\ref{eq330}),
 $$\sup\limits_{t_0-R^2 \leqslant t\leqslant t_0 } \left|B_R(x_0)\bigcap\left[u(\cdot,t)>\mu-\frac{H}{2^{s_0}}\right]\right| \leqslant \left[ 1-\sigma +\frac{1}{2}\min\{\sigma,1-\sigma\} \right]|B_R|.$$
 Employing Lemma \ref{Le32} ($s$ and $H$ are changed into $s-s_0-1$ and $\frac{H}{2^{s_0}}$ , respectively), it follows $$\left|\mathcal{Q}_R \bigcap\left[u>\mu-\frac{H}{2^{s-1}}\right]\right|\leqslant \frac{C}{\sigma\sqrt{s-s_0-1}}|\mathcal{Q}_R|,$$
 where $C$ relies on $Q,\lambda_0,\eta,$ and $\gamma(\cdot)$. Let $\theta$ be the constant in Lemma \ref{Le34} and choose $s$ large enough satisfying $$\frac{C}{\sigma\sqrt{s-s_0-1}}\leqslant\theta.$$
It implies (\ref{eq332}) from Lemma \ref{Le34}.

%---------------------------------------------------------------------------------
 \section{Proofs of main results}\label{Section 4}
%---------------------------------------------------------------------------------

We first prove an oscillation estimate for the weak solution to (\ref{eq11}).
\begin{lemma}\label{Le41}
Let $u\in{V_{2}^{1,0}({\mathcal{Q}}_{T})}$ be the weak solution to (\ref{eq11}) with $(C1),(C2)$ and $(C3)$,   $${{\mathcal{Q}}_{R_0}}={B_{R_0}(x_0)}\times{({t_0}-R^2,{t_0}]}\subset{{\mathcal{Q}}_{T}},0<R_0\leqslant1,$$ then for any $R\in(0.R_0],$ there exists $\beta,0<\beta\leqslant{1-\frac{Q+2}{\eta}},$ such that
\begin{equation}\label{eq41}
 \mathop{osc}\limits_{\mathcal{Q}_R}u\leqslant C \left(\frac{R}{R_0}\right)^{\beta}\left[\mathop{osc}\limits_{\mathcal{Q}_{R_0}}u +(M+F_0)R_0|\mathcal{Q}_{R_0}|^{-\frac{1}{\eta}}\right],
\end{equation} where $\mathop{osc}\limits_{\mathcal{Q}_R}u = ess\sup\limits_{\mathcal{Q}_R}u-ess\inf\limits_{\mathcal{Q}_R}u,C\geqslant1$ relies on $Q,\eta, \Lambda,\delta$ and $$F_0={\sum\limits_i}\|f^i\|_{{L^p}(\mathcal{Q}_T)} + \|f\|_{L^{\frac{p(Q+2)}{Q+2+p}(\mathcal{Q}_T)}}.$$
\end{lemma}
\textbf{Proof.} Denote $\nu(R)=ess\sup\limits_{\mathcal{Q}_R}u,\tilde{\nu}(R)=ess\inf\limits_{\mathcal{Q}_R}u$ and  $\omega(R)=\nu(R)-\tilde{\nu}(R)$, then $$ \mathop{osc}\limits_{\mathcal{Q}_R}u=\omega(R) $$ and one of the following two inequalities holds:
\begin{equation}\label{eq42}
\left|B_{R/2}(x_0)\bigcap\left[u\left(\cdot,t_0-\left(\frac{R}{2}\right)^2\right)<\tilde{\nu}(R) + \frac{1}{2}\omega(R)\right]\right|\leqslant\frac{1}{2}|B_{R/2}|,
\end{equation}
\begin{equation}\label{eq43}
\left|B_{R/2}(x_0)\bigcap\left[u\left(\cdot,t_0-\left(\frac{R}{2}\right)^2\right)>\nu(R) - \frac{1}{2}\omega(R)\right]\right|\leqslant\frac{1}{2}|B_{R/2}|.
\end{equation}

If (\ref{eq42}) is valid, then Lemma \ref{Le36} implies that for $\frac{\omega(R)}{2}\leqslant H:=\frac{\mu-\tilde{\mu}}{2}\leqslant\delta M, $ there exists $s_1=s_1(1/2)\geqslant 1,$ such that one of the following two inequalities holds:
\begin{equation}\label{eq44}
\frac{\omega(R)}{2}\leqslant2^{s_1}(M+F_0)R|\mathcal{Q}_{R}|^{-\frac{1}{\eta}},
\end{equation}
\begin{equation}\label{eq45}
ess\inf\limits_{\mathcal{Q}_{R/4}}u\geqslant\tilde{\nu}(R)+\frac{\omega(R)}{2^{{s_1}+2}}.
\end{equation}
It sees that by (\ref{eq44}),
\begin{equation}\label{eq46}
\omega(R/4)\leqslant\omega(R)\leqslant2^{{s_1}+1}(M+F_0)R|\mathcal{Q}_{R}|^{-\frac{1}{\eta}},
\end{equation}
and by (\ref{eq45}),
\begin{align}\label{eq47}
\omega(R/4) &=ess\sup\limits_{\mathcal{Q}_{R/4}}u-ess\inf\limits_{\mathcal{Q}_{R/4}}u \nonumber \\ &\leqslant ess\sup\limits_{\mathcal{Q}_{R}}u -\tilde{\nu}(R)-\frac{\omega(R)}{2^{{s_1}+2}} \nonumber \\ &=\omega(R)-\frac{\omega(R)}{2^{{s_1}+2}} \nonumber \\ &=\omega(R)(1-2^{{s_1}+2}).
\end{align}

If (\ref{eq43}) is valid, then by Lemma \ref{Le36} there exists $s_2=s_2(1/2)\geqslant 1,$ such that one of the following two inequalities holds:
\begin{equation}\label{eq48}
\frac{\omega(R)}{2}\leqslant2^{s_2}(M+F_0)R|\mathcal{Q}_{R}|^{-\frac{1}{\eta}},
\end{equation}
\begin{equation}\label{eq49}
ess\sup\limits_{\mathcal{Q}_{R/4}}u\leqslant\tilde{\nu}(R)-\frac{\omega(R)}{2^{{s_2}+2}}.
\end{equation}
It shows that from (\ref{eq48}),
\begin{equation}\label{eq410}
\omega(R/4)\leqslant\omega(R)\leqslant2^{{s_2}+1}(M+F_0)R|\mathcal{Q}_{R}|^{-\frac{1}{\eta}},
\end{equation}
and from (\ref{eq49}),
\begin{align}\label{eq411}
\omega(R/4) &=ess\sup\limits_{\mathcal{Q}_{R/4}}u-ess\inf\limits_{\mathcal{Q}_{R/4}}u \nonumber \\ &\leqslant \nu(R)-\frac{\omega(R)}{2^{{s_2}+2}} - ess\inf\limits_{\mathcal{Q}_{R/4}}u \nonumber \\ &=\nu(R)-\tilde{\nu}(R)-\frac{\omega(R)}{2^{{s_2}+2}} \nonumber \\ &=\omega(R)-\frac{\omega(R)}{2^{{s_2}+2}} \nonumber \\ &=\omega(R)(1-2^{{s_1}+2}).
\end{align}
Let us take $s_0=\max\{s_1,s_2\}$, we derive by (\ref{eq46}) and (\ref{eq410}) that
\begin{equation}\label{eq412}
\omega(R/4)\leqslant\omega(R)\leqslant2^{{s_0}+1}(M+F_0)R|\mathcal{Q}_{R}|^{-\frac{1}{\eta}},
\end{equation}
and by (\ref{eq47}) and (\ref{eq411}) that
\begin{equation}\label{eq413}
\omega(R/4)=ess\sup\limits_{\mathcal{Q}_{R/4}}u-ess\inf\limits_{\mathcal{Q}_{R/4}}u \leqslant \omega(R)(1-2^{{s_1}+2}).
\end{equation}

Combining (\ref{eq412}) and (\ref{eq413}), it follows that for $ R\in(0,R_0],0<R_0\leqslant 1, $
\begin{equation}\label{eq414}
\omega(R/4)=\mathop{osc}\limits_{\mathcal{Q}_{R/4}}u \leqslant \omega(R)(1-2^{{s_0}+2})+2^{{s_0}+1}(M+F_0)R|\mathcal{Q}_{R}|^{-\frac{1}{\eta}}.
\end{equation}
In terms of $|B_{{R_0}/4}|\geqslant C_1^{-1}\left(\frac{1}{4}\right)^Q|B_{R_0}|,$ we find $$\frac{R_0}{4}|\mathcal{Q}_{{R_0}/4}|^{-\frac{1}{\eta}}=\frac{R_0}{4} \left[\left(\frac{R_0}{4}\right)^2|B_{{R_0}/4}|\right]^{-\frac{1}{\eta}}\leqslant C_1^{\frac{1}{\eta}}\left(\frac{1}{4}\right)^{1-\frac{Q+2}{\eta}}R_0|\mathcal{Q}_{R_0}|^{-\frac{1}{\eta}}.$$ Replacing $R$ in (\ref{eq414}) by $R_0$ and denoting $\vartheta=1-2^{-({s_0}+1)}$ and $\upsilon=\left(\frac{1}{4}\right)^{1-\frac{Q+2}{\eta}}$, it implies $0<\vartheta,\upsilon<1$, and from (\ref{eq414}) and the above estimate that
\begin{align*}
\omega({R_0}/{4^2}) &\leqslant \vartheta\omega({R_0}/4)+2^{{s_0}+1}(M+F_0)\frac{R_0}{4}|\mathcal{Q}_{{R_0}/4}|^{-\frac{1}{\eta}} \\ &\leqslant\vartheta\left(\vartheta\omega(R_0)+2^{{s_0}+1}(M+F_0)R_0|\mathcal{Q}_{R_0}|^{-\frac{1}{\eta}}\right) + 2^{{s_0}+1}(M+F_0)C_1^{\frac{1}{\eta}}\upsilon R_0|\mathcal{Q}_{R_0}|^{-\frac{1}{\eta}} \\ &\leqslant \vartheta^2\omega(R_0)+(\vartheta+\upsilon)2^{{s_0}+1}(M+F_0)C_1^{\frac{1}{\eta}} {R_0}|\mathcal{Q}_{R_0}|^{-\frac{1}{\eta}}.
\end{align*}
Generally, for $\ell=\log_4{\frac{R_0}{R}}$, we have
\begin{equation}\label{eq415}
\omega({R_0}/{4^\ell})\leqslant\vartheta^{\ell}\omega(R_0) + \sum\limits_{\chi=0}^{\ell-1}\vartheta^{\chi}\upsilon^{\ell-1-\chi}2^{{s_0}+1}C_1^{\frac{1}{\eta}}(M+F_0) {R_0}|\mathcal{Q}_{R_0}|^{-\frac{1}{\eta}}.
\end{equation}Let us discuss two cases: $\vartheta\geqslant\upsilon$ and $\vartheta<\upsilon$.

(1) If $\vartheta\geqslant\upsilon$, it gives $\frac{1}{\vartheta}\leqslant\frac{1}{\upsilon}=4^{1-\frac{Q+2}{\eta}}$ and then $\log_4{\frac{1}{\vartheta}}\leqslant 1-\frac{Q+2}{\eta}$. Choosing $\beta_1=\log_4{\frac{1}{\vartheta}},$ we derive from (\ref{eq415}) that
\begin{align*}
\omega(R) &\leqslant\vartheta^{\ell}\omega(R_0)+\frac{1-\vartheta^{\ell}}{1-\vartheta}2^{{s_0}+1}C_1^{\frac{1}{\eta}}(M+F_0) {R_0}|\mathcal{Q}_{R_0}|^{-\frac{1}{\eta}} \\ &\leqslant \vartheta^{\ell}\left(\omega(R_0)+\frac{(\vartheta^{-\ell}-1)2^{{s_0}+1}C_1^{\frac{1}{\eta}}}{1-\vartheta}(M+F_0) {R_0}|\mathcal{Q}_{R_0}|^{-\frac{1}{\eta}}\right) \\ &\leqslant C\left(\frac{R}{R_0}\right)^{\beta_1}\left[\omega(R_0)+(M+F_0){R_0}|\mathcal{Q}_{R_0}|^{-\frac{1}{\eta}}\right],
\end{align*}
where we used $\vartheta^{\ell}=\left(\frac{R}{R_0}\right)^{\beta_1}$ and $C=\max\left\{1,\frac{(\vartheta^{-\ell}-1)2^{{s_0}+1}C_1^{\frac{1}{\eta}}}{1-\vartheta}\right\}.$

(2)  If $\vartheta<\upsilon$, then $\frac{1}{\vartheta}>\frac{1}{\upsilon}=4^{1-\frac{Q+2}{\eta}}$ and $\log_4{\frac{1}{\vartheta}}> 1-\frac{Q+2}{\eta}$. Choosing $\beta_2=1-\frac{Q+2}{\eta},$ it leads to from (\ref{eq415}) that
\begin{align*}
\omega(R) &\leqslant\upsilon^{\ell}\omega(R_0)+\frac{\upsilon^{\ell}-1}{\upsilon-1}2^{{s_0}+1}C_1^{\frac{1}{\eta}}(M+F_0) {R_0}|\mathcal{Q}_{R_0}|^{-\frac{1}{\eta}} \\ &\leqslant \upsilon^{\ell}\left(\omega(R_0)+\frac{(1-\upsilon^{-\ell})2^{{s_0}+1}C_1^{\frac{1}{\eta}}}{\upsilon-1}(M+F_0) {R_0}|\mathcal{Q}_{R_0}|^{-\frac{1}{\eta}}\right) \\ &\leqslant C\left(\frac{R}{R_0}\right)^{\beta_2}\left[\omega(R_0)+(M+F_0){R_0}|\mathcal{Q}_{R_0}|^{-\frac{1}{\eta}}\right],
\end{align*}
where $\upsilon^{\ell}=\left(\frac{R}{R_0}\right)^{1-\frac{Q+2}{\eta}}$ and $C=\max\left\{1,\frac{(1-\upsilon^{-\ell})2^{{s_0}+1}C_1^{\frac{1}{\eta}}}{\upsilon-1}\right\}.$

Select $\beta=\min\{\beta_1,\beta_2\},$ we prove (\ref{eq41}) by combining estimates in cases (1) and (2).

The following is an isomorphism between $\mathcal{L}^{p',\lambda}(\mathcal{Q})$ and $C^\alpha(\mathcal{Q}).$

\begin{lemma}\label{Le42} Let $\mathcal{Q}\subset\subset\mathcal{Q}_T, $ we have for $0<\lambda\leqslant p',p'\geqslant 1,$ \[\mathcal{L}^{p',\lambda}(\mathcal{Q})\cong C^\alpha(\mathcal{Q}),\] where $\alpha=\frac{\lambda}{p'}$.
\end{lemma}
\textbf{Proof.} Suppose $u\in{C^{\alpha}(\mathcal{Q})}.$ For any $Z=(x,t)\in\mathcal{Q},$ denote
\[\hat{\mathcal{Q}}_R(Z)=\left(B_R(x_0)\times (t,t+R^2]\right)\bigcap{\mathcal{Q}}, ~for~ 0<R\leqslant d={\rm{diam}} {\kern 1pt}{\mathcal{Q}},\] we have for $Y_1=(y_1,s_1)$ and $Y_2=(y_2,s_2)\in{\hat{\mathcal{Q}}_R(Z)}$,
\begin{align*}
\left|u(Y_1)-u_{\hat{\mathcal{Q}}_R(Z)}\right| &\leqslant \frac{1}{|\hat{\mathcal{Q}}_R(Z)|}\int_{\hat{\mathcal{Q}}_R(Z)}\left|u(Y_1)-u(Y_2)\right| dY_2 \\ &\leqslant \frac{[u]_\alpha}{|\hat{\mathcal{Q}}_R(Z)|}\int_{\hat{\mathcal{Q}}_R(Z)}d_{\mathcal{P}}^{\alpha}(Y_1,Y_2)dY_2 \\ &\leqslant C[u]_\alpha R^\alpha
\end{align*}
and $$\frac{R^{-\lambda}}{|\hat{\mathcal{Q}}_R(Z)|}\int_{\hat{\mathcal{Q}}_R(Z)} \left|u(Y_1)-u_{\hat{\mathcal{Q}}_R(Z)}\right|^{p'}dY_1\leqslant CR^{\alpha p'-\lambda}[u]_{\alpha}^{p'}.$$
Noting $\alpha=\frac{\lambda}{p'},$ it yields $$\left\{\sup\limits_{Z\in\mathcal{Q},0<R\leqslant d}\frac{R^{-\lambda}}{|\hat{\mathcal{Q}}_R(Z)|}\int_{\hat{\mathcal{Q}}_R(Z)} \left|u(Y_1)-u_{\hat{\mathcal{Q}}_R(Z)}\right|^{p'}dY_1\right\}^{\frac{1}{p'}}\leqslant CR^{\alpha-\frac{\lambda}{p'}}[u]_\alpha=C[u]_\alpha$$ and
\begin{equation}\label{eq416}
C^\alpha(\mathcal{Q})\subset\mathcal{L}^{p',\lambda}(\mathcal{Q}).
\end{equation}

On the contrary, if $u\in\mathcal{L}^{p',\lambda}(\mathcal{Q})$, we have for any $Y\in\hat{\mathcal{Q}}_R(Z),0<\rho<R\leqslant d,$
\begin{equation}\label{eq417}
\left|u_{\hat{\mathcal{Q}}_\rho(Z)}-u_{\hat{\mathcal{Q}}_R(Z)}\right|^{p'} \leqslant 2^{p'-1}\left[\left|u_{\hat{\mathcal{Q}}_\rho(Z)}-u(Y)\right|^{p'} + \left|u(Y)-u_{\hat{\mathcal{Q}}_R(Z)}\right|^{p'}\right].
\end{equation}
Integrating it over ${\hat{\mathcal{Q}}_\rho(Z)}$, it follows
\begin{align*}
\left|u_{\hat{\mathcal{Q}}_\rho(Z)}-u_{\hat{\mathcal{Q}}_R(Z)}\right|^{p'} \left|{\hat{\mathcal{Q}}_\rho(Z)}\right| &\leqslant 2^{p'-1}\left[\int_{{\hat{\mathcal{Q}}_\rho(Z)}} u_{\hat{\mathcal{Q}}_\rho(Z)}-\left|u(Y)\right|^{p'}dY + \int_{{\hat{\mathcal{Q}}_R(Z)}} \left|u(Y)-u_{\hat{\mathcal{Q}}_R(Z)}\right|^{p'}dY \right] \\ &\leqslant 2^{p'}R^\lambda\left|{\hat{\mathcal{Q}}_\rho(Z)}\right|[u]_{p',\lambda}^{p'}
\end{align*}
and $$\left|u_{\hat{\mathcal{Q}}_\rho(Z)}-u_{\hat{\mathcal{Q}}_R(Z)}\right|\leqslant CR^{\frac{\lambda}{p'}}\left|\frac{{\hat{\mathcal{Q}}_R(Z)}}{{\hat{\mathcal{Q}}_\rho(Z)}}\right|^{\frac{1}{p'}} [u]_{p',\lambda}.$$
Using $\left|\frac{{\hat{\mathcal{Q}}_R(Z)}}{{\hat{\mathcal{Q}}_\rho(Z)}}\right| \leqslant\left(\frac{R}{\rho}\right)^{Q+2}$ by (\ref{eq22}), it views that for any $k,m,k\leqslant m,$
\begin{align}\label{eq418}
\left|u_{\hat{\mathcal{Q}}_{2^{-m}R}(Z)}-u_{\hat{\mathcal{Q}}_{2^{-k}R}(Z)}\right| &\leqslant \sum\limits_{j=k+1}^{m}\left|u_{\hat{\mathcal{Q}}_{2^{-j}R}(Z)}-u_{\hat{\mathcal{Q}}_{2^{-j+1}R}(Z)}\right| \nonumber \\ &\leqslant \sum\limits_{j=k+1}^{m} C\left|{2^{-j+1}R}\right|^{\alpha}2^{\frac{Q+2}{p'}}[u]_{p',\lambda} \nonumber \\ &=C2^{\alpha+\frac{Q+2}{p'}}R^\alpha \sum\limits_{j=k+1}^{m} 2^{-j\alpha}[u]_{p',\lambda} \nonumber \\ &\leqslant CR^\alpha|2^{-k\alpha}-2^{-m\alpha}|[u]_{p',\lambda},
\end{align}
which implies that $\left\{u_{\hat{\mathcal{Q}}_{2^{-m}R}(Z)}\right\}$ is a Cauchy sequence, and its limit is denoted by $\tilde{u}(Z)$. Letting $m\to\infty$ in (\ref{eq418}), we have
\begin{equation}\label{eq419}
\left|u_{\hat{\mathcal{Q}}_{2^{-k}R}(Z)}-\tilde{u}(Z)\right|\leqslant CR^{\alpha}2^{-k\alpha}[u]_{p',\lambda}.
\end{equation}
Since $u_{\hat{\mathcal{Q}}_{2^{-k}R}(Z)}$ (for any $k$) is continuous on $Z$, it knows by (\ref{eq419}) that $\tilde{u}(Z)$ is continuous in $\mathcal{Q}$. In addition, Lebesgue's Theorem assures $$u_{\hat{\mathcal{Q}}_{2^{-k}R}(Z)}\rightarrow u(Z),~a.e. Z\in\mathcal{Q},$$
hence $u(Z)=\tilde{u}(Z)$ a.e. $Z\in\mathcal{Q}$ and $\tilde{u}(Z)$ is independent of $R$. Now we see from (\ref{eq419}) that
\begin{equation}\label{eq420}
\left|u_{\hat{\mathcal{Q}}_{R}(Z)}-\tilde{u}(Z)\right|\leqslant CR^{\alpha}[u]_{p',\lambda}.
\end{equation}

Let us prove further that $\tilde{u}(Z)$ is H\"{o}lder continuous. For any $Y_1,Y_2\in\mathcal{Q}$, we take $R=d_{\mathcal{P}}(Y_1,Y_2)$ and obtain from (\ref{eq420}) that
\begin{align}\label{eq421}
\left|\tilde{u}(Y_1)-\tilde{u}(Y_2)\right| &\leqslant \left|u_{\hat{\mathcal{Q}}_{R}(Y_1)}-\tilde{u}(Y_1)\right| + \left|u_{\hat{\mathcal{Q}}_{R}(Y_2)}-\tilde{u}(Y_2)\right| + \left|u_{\hat{\mathcal{Q}}_{R}(Y_1)}-u_{\hat{\mathcal{Q}}_{R}(Y_2)}\right| \nonumber \\ &\leqslant  CR^{\alpha}[u]_{p',\lambda} + \left|u_{\hat{\mathcal{Q}}_{R}(Y_1)}-u_{\hat{\mathcal{Q}}_{R}(Y_2)}\right|.
\end{align}
To estimate $ \left|u_{\hat{\mathcal{Q}}_{R}(Y_1)}-u_{\hat{\mathcal{Q}}_{R}(Y_2)}\right|$ in (\ref{eq421}), we note
 \[ \left|u_{\hat{\mathcal{Q}}_{R}(Y_1)}-u_{\hat{\mathcal{Q}}_{R}(Y_2)}\right|^{p'}\leqslant 2^{p'-1}\left[\left|u_{\hat{\mathcal{Q}}_R(Y_1)}-\tilde{u}(Z')\right|^{p'} + \left|u_{\hat{\mathcal{Q}}_R(Y_2)}-\tilde{u}(Z')\right|^{p'}\right], for~ any~ Z'\in\mathcal{Q},\] and integrate it over ${\hat{\mathcal{Q}}_R(Y_1)}\cap{\hat{\mathcal{Q}}_R(Y_2)}$ to gain
 \begin{align*}
 & ~\left|{\hat{\mathcal{Q}}_R(Y_1)}\cap{\hat{\mathcal{Q}}_R(Y_1)}\right| \left|u_{\hat{\mathcal{Q}}_{R}(Y_1)}-u_{\hat{\mathcal{Q}}_{R}(Y_2)}\right|^{p'} \\ &\leqslant 2^{p'-1}\left[\int_{{\hat{\mathcal{Q}}_R(Y_1)}\cap{\hat{\mathcal{Q}}_R(Y_1)}} \left|u_{\hat{\mathcal{Q}}_R(Y_1)}-\tilde{u}(Z')\right|^{p'}dZ' + \int_{{\hat{\mathcal{Q}}_R(Y_1)}\cap{\hat{\mathcal{Q}}_R(Y_1)}} \left|u_{\hat{\mathcal{Q}}_R(Y_2)}-\tilde{u}(Z')\right|^{p'}dZ' \right] \\ &\leqslant CR^\lambda \left|{\hat{\mathcal{Q}}_R(Y_1)}\cap{\hat{\mathcal{Q}}_R(Y_1)}\right|[u]_{p',\lambda}^{p'},
 \end{align*}
 which shows $$\left|u_{\hat{\mathcal{Q}}_{R}(Y_1)}-u_{\hat{\mathcal{Q}}_{R}(Y_2)}\right| \leqslant CR^{\alpha}[u]_{p',\lambda}.$$
 Substituting it into (\ref{eq421}), it happens $$\left|\tilde{u}(Y_1)-\tilde{u}(Y_2)\right|\leqslant CR^{\alpha}[u]_{p',\lambda}$$ and noting $R=d_{\mathcal{P}}(Y_1,Y_2),$ we arrive at
 \[\mathcal{L}^{p',\lambda}(\mathcal{Q})\subset C^\alpha(\mathcal{Q}). \]

\textbf{Proof of Theorem \ref{Th11}.} For any $Z_0=(x_0,t_0)\in\mathcal{Q}$, denote $$\mathcal{Q}_R(Z_0)= B_R(x_0)\times ({t_0}-{R^2},t_0],$$ By employing Lemma \ref{Le41}, we derive for $R\in(0,\bar{d}]$, $\bar{d}=\min\{1,d_{\mathcal{P}}(\mathcal{Q},\partial_{\mathcal{P}}\mathcal{Q}_T)\}$,
\begin{align}\label{eq422}
\mathop{osc}\limits_{\mathcal{Q}_R(Z_0)} u &\leqslant C\left(\frac{R}{\bar{d}}\right)^{\beta}\left[\mathop{osc}\limits_{\mathcal{Q}_T} u +(M+F_0){\bar{d}}|\mathcal{Q}_T|^{-\frac{1}{\eta}}\right] \nonumber \\ &\leqslant C\left(\frac{R}{\bar{d}}\right)^{\beta}\left[ M \left(1+{\bar{d}}|\mathcal{Q}_T|^{-\frac{1}{\eta}}\right) + F_0{\bar{d}}|\mathcal{Q}_T|^{-\frac{1}{\eta}}\right] \nonumber \\ &\leqslant C\left(\frac{R}{\bar{d}}\right)^{\beta}\left[ M + F_0{\bar{d}}|\mathcal{Q}_T|^{-\frac{1}{\eta}}\right],
\end{align}
where $ 1+{\bar{d}}|\mathcal{Q}_T|^{-\frac{1}{\eta}}$ is a constant, and then $$\frac{R^{-\beta}}{|\mathcal{Q}_R(Z_0)|}\int_{\mathcal{Q}_R(Z_0)}\left|u(Z)-u_{\mathcal{Q}_R(Z_0)}\right|dZ \leqslant R^{-\beta}\mathop{osc}\limits_{\mathcal{Q}_R(Z_0)} u \leqslant C{\bar{d}}^\beta \left(M + F_0{\bar{d}}|\mathcal{Q}_T|^{-\frac{1}{\eta}}\right),$$ which implies $u\in\mathcal{L}_{loc}^{1,\beta}(\mathcal{Q}_T;d_{\mathcal{P}})$ and so $u\in C_{loc}^\beta(\mathcal{Q}_T;d_{\mathcal{P}})$ by Lemma \ref{Le42}. Now let us estimate $[u]_{\beta;\mathcal{Q}}$.

For any $Z_0=(x_0,t_0),Z_1=(x_1,t_1)\in\mathcal{Q}$, without losing of generality, suppose $t_0\geqslant t_1$. If $d_{\mathcal{P}}(Z_0,Z_1)\leqslant\bar{d}$, then $Z_1\in\mathcal{Q}_R(Z_0)\subset\mathcal{Q}_T,R=d_{\mathcal{P}}(Z_0,Z_1)$, and from (\ref{eq422}),
\begin{equation}\label{eq423}
\left|u(Z_0)-u(Z_1)\right|\leqslant\mathop{osc}\limits_{\mathcal{Q}_R(Z_0)}u \leqslant C\left(\frac{d_{\mathcal{P}}(Z_0,Z_1)}{\bar{d}}\right)^\beta\left(M + F_0{\bar{d}}|\mathcal{Q}_T|^{-\frac{1}{\eta}}\right);
\end{equation}
if $d_{\mathcal{P}}(Z_0,Z_1)>\bar{d}$, then
\begin{equation}\label{eq424}
\left|u(Z_0)-u(Z_1)\right|\leqslant\frac{2M}{{\bar{d}}^\beta}\left[d_{\mathcal{P}}(Z_0,Z_1)\right]^\beta.
\end{equation}
 It follows (\ref{eq14}) by combining (\ref{eq423}) and (\ref{eq424}).

To prove Theorem \ref{Th12}, we need an extension property of positivity for functions in the De Giorgi class.
\begin{lemma}\label{Le43} Let $u\in{DG^-}(\mathcal{Q}_T)$ and $u\geqslant0$ in $$ \mathcal{Q}_R^a=B_R(x_0)\times (t_0,t_0+aR^2)\subset\mathcal{Q}_T, 0<R\leqslant 1.$$ For $\epsilon\in(0,1)$, if
\begin{equation}\label{eq425}
ess\inf\limits_{B_{\epsilon R}}u(x,t_0)\geqslant k\geqslant 0,
\end{equation} then there exist $R_0>0$ and a positive integer $s\geqslant1$ relying on $Q,\lambda_0,\eta$ and $\gamma(\cdot)$, such that for $B_R=B_R(x_0),0<R\leqslant R_0$, either
\begin{equation}\label{eq426}
k\leqslant 2^{s+2}{F_0}R|\mathcal{Q}_{R}|^{-\frac{1}{\eta}},
\end{equation} or
\begin{equation}\label{eq427}
ess\inf\limits_{B_{R/4}} u(x,t+aR^2)\geqslant \epsilon^s(k-2^{s+2}{F_0}R|\mathcal{Q}_{R}|^{-\frac{1}{\eta}}).
\end{equation}
\end{lemma}
\textbf{Proof.} Let $\epsilon\leqslant1/8$, we have from (\ref{eq425}) and (\ref{eq22}) that $$\left|B_{4\epsilon R}\bigcap [u(\cdot,t_0)>k]\right|\geqslant |B_{\epsilon R}|$$ and then $$\left|B_{4\epsilon R}\bigcap [u(\cdot,t_0)<k]\right|\leqslant|B_{4\epsilon R}|-\left|B_{4\epsilon R}\bigcap [u(\cdot,t_0)>k]\right|\leqslant|B_{4\epsilon R}|-|B_{\epsilon R}|\leqslant (1-C_1^{-1}4^{-Q})|B_{4\epsilon R}|. $$ As a result of Lemma \ref{Le36}, there exists $s>1$, such that either
\begin{equation}\label{eq428}
k\leqslant 2^s(k+F_0)(4\epsilon R)|\mathcal{Q}_{4\epsilon R}|^{-\frac{1}{\eta}},
\end{equation} ($M$ in (\ref{eq331}) is changed to $k$; actually, it is suitable from the proofs of lemmas in Section 3), or
\begin{equation}\label{eq429}
ess\inf\limits_{\mathcal{Q}_{2\epsilon R}} u\geqslant \frac{k}{2^s},
\end{equation} where we used $H=k$ which is followed from $H:=k-\inf u$ and $u\geqslant 0$.

If (\ref{eq428}) holds, we choose $R_0$ satisfying $$2^{s+1}R_0|\mathcal{Q}_{R_0}|^{-\frac{1}{\eta}}=\frac{1}{2}$$ and for $0<R\leqslant R_0$, it yields by (\ref{eq22}) that $|\mathcal{Q}_{4\epsilon R}|^{-\frac{1}{\eta}}\leqslant C_1^{\frac{1}{\eta}}\frac{R_0}{4\epsilon R}|\mathcal{Q}_{R_0}|^{-\frac{1}{\eta}}$ and $|\mathcal{Q}_{4\epsilon R}|^{-\frac{1}{\eta}}\leqslant C_1^{\frac{1}{\eta}}\frac{R}{4\epsilon R}|\mathcal{Q}_{R}|^{-\frac{1}{\eta}}.$ Noting $C_1^{\frac{1}{\eta}}=2^{\frac{Q}{\eta}}<2,$ it implies from (\ref{eq428}) that
\begin{align*}
k &\leqslant 2^sk(4\epsilon R)|\mathcal{Q}_{4\epsilon R}|^{-\frac{1}{\eta}} + 2^s{F_0}(4\epsilon R)|\mathcal{Q}_{4\epsilon R}|^{-\frac{1}{\eta}} \\ &\leqslant {2^s}k{C_1^{\frac{1}{\eta}}}R_0|\mathcal{Q}_{R_0}|^{-\frac{1}{\eta}} + {2^s}C_1^{\frac{1}{\eta}}{F_0}R|\mathcal{Q}_{R}|^{-\frac{1}{\eta}} \\ &\leqslant \frac{1}{2} k + {2^{s+1}}{F_0}R|\mathcal{Q}_{R}|^{-\frac{1}{\eta}}.
\end{align*} This proves (\ref{eq426}).

If (\ref{eq429}) holds, choose a suitable $\epsilon>0$ such that $N=\frac{\log{\epsilon^{-1}}}{\log2}-2$ is an integer, and arrive at by employing Lemma \ref{Le36} with $N$ times that
\begin{align*}
ess\inf\limits_{B_{R/4}} u(x,t_0 +aR^2) &=ess\inf\limits_{B_{{2^N}\epsilon R}} u(x,t_0 +aR^2) \\ &\geqslant ess\inf\limits_{\mathcal{Q}_{{2^N}\epsilon R}} u \geqslant \frac{k}{2^{Ns}} \geqslant \epsilon^s k \\ &\geqslant \epsilon^s\left(k-2^{s+2}{F_0}R|\mathcal{Q}_{R}|^{-\frac{1}{\eta}}\right)
\end{align*}
which is (\ref{eq427}).

\textbf{Proof of Theorem \ref{Th12}.} For simplicity, we assume $a'=2$ and denote $$\mathcal{Q}_r=B_r(x_0)\times(t_0+R^2-r^2,t_0+R^2+r^2).$$ Consider two functions
\[ m(r)=u(x_0,t_0+R^2)\left(1-\frac{r}{R}\right)^{-s} ~and~\mu(r)=\max\limits_{\mathcal{Q}_r}u(x,t), \]
where $s$ is determined in Lemma \ref{Le43}, and let $r_0 (r_0<R)$ be the largest root to the equation $m(r)=\mu(r).$ Since $m(r)\to\infty$ as $r\to R-0,$ and $u$ is continuous and bounded in $\mathcal{Q}_R^2$, it sees that $r_0$ is well defined, $m(r)>\mu(r)$ for $r_0<r\leqslant R,$ and there exists $(x_1,t_1)\in\mathcal{Q}_{r_0}=B_{r_0}(x_0)\times(t_0+R^2-r_0^2,t_0+R^2+r_0^2)$, such that  $$u(x_1,t_1)=m(r_0)=\mu(r_0).$$ Now introduce
\[ \mathcal{Q}=\left\{(x,t)\mid d(x,x_1)\leqslant\frac{R-r_0}{2},t_1-\frac{(R-r_0)^2}{4}<t\leqslant t_1\right\}, \]
i.e. $\mathcal{Q}=B_{(R-r_0)/2}(x_1)\times\left( t_1-\frac{(R-r_0)^2}{4},t_1 \right].$ For any $(x,t)\in\mathcal{Q}$, \[d(x,x_0)\leqslant d(x,x_1)+d(x_0,x_1)\leqslant\frac{R-r_0}{2}+r_0=\frac{R+r_0}{2}\] and
\begin{align}\label{eq430}
t_0+R^2-\frac{(R+r_0)^2}{4} &<t_0+R^2-r_0^2-\frac{(R-r_0)^2}{4} \nonumber \\ &<t_1-\frac{(R-r_0)^2}{4}<t_1 \nonumber \\ &<t_0+R^2+r_0^2 <t_0+R^2+\frac{(R+r_0)^2}{4},
\end{align} where $r_0^2<r_0^2+\frac{(R-r_0)^2}{4}<\frac{(R+r_0)^2}{4},$ then we have $\mathcal{Q}\subset\mathcal{Q}_{(R+r_0)/2}.$ Noting the meaning of $r_0$, it follows
\begin{align}\label{eq431}
\sup\limits_{\mathcal{Q}}u(x,t) &\leqslant\mu\left(\frac{R+r_0}{2}\right)\leqslant m \left(\frac{R+r_0}{2}\right) \nonumber \\ &=u(x_0,t_0+R^2)\cdot2^s\left(1-\frac{r_0}{R}\right)^{-s} \nonumber \\ &=2^s m(r_0) \nonumber \\ &=2^s \mu(r_0).
\end{align} From Lemma \ref{Le41}, (\ref{eq22}), (\ref{eq431}) and $\frac{R-r_0}{2}<R $, we obtain
\begin{align*}
& |u(x_1,t_1)-\inf\limits_{d(x,x_1)\leqslant{\epsilon(R-r_0)}}u(x,t_1)|\leqslant \mathop{osc}\limits_{B_{\epsilon(R-r_0)}} u(x,t_1) \\ &\leqslant C2^{\beta}{\epsilon}^{\beta}\left[\mathop{osc}\limits_{B_{(R-r_0)/2}} u(x,t_1) + F_0\frac{R-r_0}{2}|\mathcal{Q}_{\frac{R-r_0}{2}}|^{-\frac{1}{\eta}}\right] \\ &\leqslant C{\epsilon}^{\beta}\left[\sup\limits_{\mathcal{Q}}u + F_0(R-r_0)|\mathcal{Q}_{R-r_0}|^{-\frac{1}{\eta}}\right] \\ &\leqslant C{\epsilon}^{\beta}\left[{2^s}\mu(r_0)+{F_0}R|\mathcal{Q}_{R}|^{-\frac{1}{\eta}}\right]
\end{align*} and
\begin{align*}
\inf\limits_{d(x,x_1)\leqslant{\epsilon(R-r_0)}}u(x,t_1) &\geqslant u(x_1,t_1)-C\epsilon^\beta\left[{2^s}\mu(r_0)+{F_0}R|\mathcal{Q}_{R}|^{-\frac{1}{\eta}}\right] \\ &=\mu(r_0)-C\epsilon^\beta\left[{2^s}\mu(r_0)+{F_0}R|\mathcal{Q}_{R}|^{-\frac{1}{\eta}}\right].
\end{align*}
Choosing $\epsilon$ so small that $C\epsilon^\beta<1$ and $C\epsilon^\beta2^s=\frac{1}{2},$ it yields
\begin{equation}\label{eq432}
\inf\limits_{d(x,x_1)\leqslant{\epsilon(R-r_0)}}u(x,t_1)\geqslant \frac{1}{2}\mu(r_0)-{F_0}R|\mathcal{Q}_{R}|^{-\frac{1}{\eta}}.
\end{equation}
For $x\in{B_{R/2}(x_0)}$, we have $$d(x,x_1)\leqslant d(x,x_0)+d(x_0,x_1)\leqslant\frac{R}{2}+r_0\leqslant\frac{3R}{2}$$ and so $$B_{R/2}(x_0)\subset{B_{3R/2}(x_1)}.$$
Employing it and (\ref{eq425}) in Lemma \ref{Le43} on the domain $B_{6R}(x_1)\times(t_1,t_0+2R^2]$, it shows
\begin{align*}
\inf\limits_{d(x,x_0)\leqslant{R/2}}u(x,t_0+2R^2) &\geqslant \inf\limits_{d(x,x_1)\leqslant{3R/2}}u(x,t_0+2R^2) \\ &\geqslant \left[\frac{\epsilon(R-r_0)}{6R}\right]^s\left(\frac{1}{2}\mu(r_0) - (2^{s+2}+1){F_0}R|\mathcal{Q}_{R}|^{-\frac{1}{\eta}}\right) \\ &\geqslant \frac{1}{2}\left[\frac{\epsilon(R-r_0)}{6R}\right]^s\cdot u(x_0,t_0+R^2)\left(1-\frac{r}{R}\right)^{-s} - \left[\frac{\epsilon(R-r_0)}{6R}\right]^s 2\cdot 6^s{F_0}R|\mathcal{Q}_{R}|^{-\frac{1}{\eta}} \\ &\geqslant \frac{1}{2}\left(\frac{\epsilon}{6}\right)^s u(x_0,t_0+R^2)-2\left[\frac{\epsilon(R-r_0)}{R}\right]^s \cdot {F_0}R|\mathcal{Q}_{R}|^{-\frac{1}{\eta}} \\ &=\frac{1}{2}\left(\frac{\epsilon}{6}\right)^s u(x_0,t_0+R^2)-2\left(\frac{6}{\epsilon}\right)^s \cdot {F_0}R|\mathcal{Q}_{R}|^{-\frac{1}{\eta}},
\end{align*}
where we used $2^{s+2}+1<2\cdot6^s$ and $\frac{\epsilon(R-r_0)}{R}<\frac{6}{\epsilon}$. This proves (\ref{eq15}).

%%appendix------------------------------------------------------------------------

%%-------------------------the followings are the proofs of e102, e103

%% bibliography--------------------------------------------------------------------

\end{document}